\documentclass[11pt,reqno]{amsart}
\usepackage{amsmath,amssymb,mathrsfs}
\usepackage{graphicx,cite,times,color}
\usepackage{cases}

\usepackage{lipsum}
\usepackage{amsfonts}
\usepackage{graphicx}
\usepackage{epstopdf}
\usepackage{algorithmic}
\usepackage{amsmath,amssymb}
\usepackage{subfigure}
\usepackage{color}
\usepackage{dcolumn}
\usepackage{bm}
\usepackage{multirow}
\usepackage{fancyhdr}
\usepackage{latexsym}
\usepackage{amscd}
\usepackage{cases}
\usepackage{algorithmic}
\usepackage{algorithm}
\usepackage{setspace}
\usepackage{booktabs}
\usepackage{enumerate}

\numberwithin{equation}{section}
\begin{document}
	\title[Sparse Representation of Gaussian Molecular Surface]{Sparse Representation of Gaussian Molecular Surface}
	
	\author{Sheng Gui}
	\address{State Key Laboratory of Scientific and Engineering Computing, National Center for Mathematics and Interdisciplinary Sciences, Academy of Mathematics and Systems Science, Chinese Academy of Sciences, Beijing 100190, China, and School of Mathematical Sciences, University of Chinese Academy of Sciences, Beijing 100049, China.}
	\email{shenggui@lsec.cc.ac.cn}
	
	\author{Minxin Chen*}
	\address{Corresponding author: Department of Mathematics, Soochow University, Suzhou 215006, China.}
	\email{chenminxin@suda.edu.cn}
	
	\author{Benzhuo Lu*}
	\address{Corresponding author: State Key Laboratory of Scientific and Engineering Computing, National Center for Mathematics and Interdisciplinary Sciences, Academy of Mathematics and Systems Science, Chinese Academy of Sciences, Beijing 100190, China, and School of Mathematical Sciences, University of Chinese Academy of Sciences, Beijing 100049, China.}
	\email{bzlu@lsec.cc.ac.cn}
	
	\thanks{This work was supported by the National Key Research and Development Program of China (grant 2016YFB0201304), Science Challenge Project (TZ2016003) and the China NSF (NSFC 21573274, NSFC 11771435), and NSF of Jiangsu Province (BK20130278).}
	
	\keywords{Gaussian surface; Molecular surface; Rotational ellipsoid feature; Radial basis function; Sparse representation; Shape comparisons}
	
	\begin{abstract}
		In this paper, we propose a model and algorithm for sparse representing Gaussian molecular surface. The original Gaussian molecular surface is approximated by a relatively small number of radial basis functions (RBFs) with rotational ellipsoid feature. The sparsity of the RBF representation is achieved by solving a nonlinear $L_1$ optimization problem. Experimental results demonstrate that the original Gaussian molecular surface is able to be represented with good accuracy by much fewer RBFs using our $L_1$ model and algorithm. The sparse representation of Gaussian molecular surface is useful in various applications, such as molecular structure alignment, calculating molecular areas and volumes, and the method in principle can be applied to sparse representation of general shapes and coarse-grained molecular modeling.

	\end{abstract}
	
	\maketitle
	
\section{Introduction}
The development of molecular surface representation in the fields of computational biology and computer-aided drug design has been a vital issue of considerable interest for many years, such as shape-based docking problems \cite{mcgann2003gaussian}, molecular shape comparisons \cite{grant1996a}, calculating SAS areas \cite{weiser1999optimization}, coarse-grained molecular dynamics \cite{Rudd2000coarse}, and the generalized Born models \cite{yu2006what}, etc. In implicit-solvent modeling (e.g., see a review in \cite{Lu2008Recent}), the molecular surface is also a prerequisite for using the boundary element method (BEM) and the finite element method (FEM).  Also the molecular shape is fundamental to molecular recognition events such as a drug binding to a biological receptor in the computer-aided drug design field \cite{Nicholls2010Molecular}. The visualisation of molecules is another field where surface representations of the molecular surface are of particular importance. Visually inspecting molecules can lead to biophysical development, since visualising 3D molecular complexes has become a common practice in life sciences \cite{Gui2018}. So far, due to the highly complex and irregular shape of the molecular surface, the efficient representation of the molecular surface for large real biomolecule with high quality remains a challenging problem \cite{Chen2013Advances}.

\begin{figure}[h]
  \centering
  \subfigure[]{
  \label{Gaussian:PQR}
  \includegraphics[scale = 0.2]{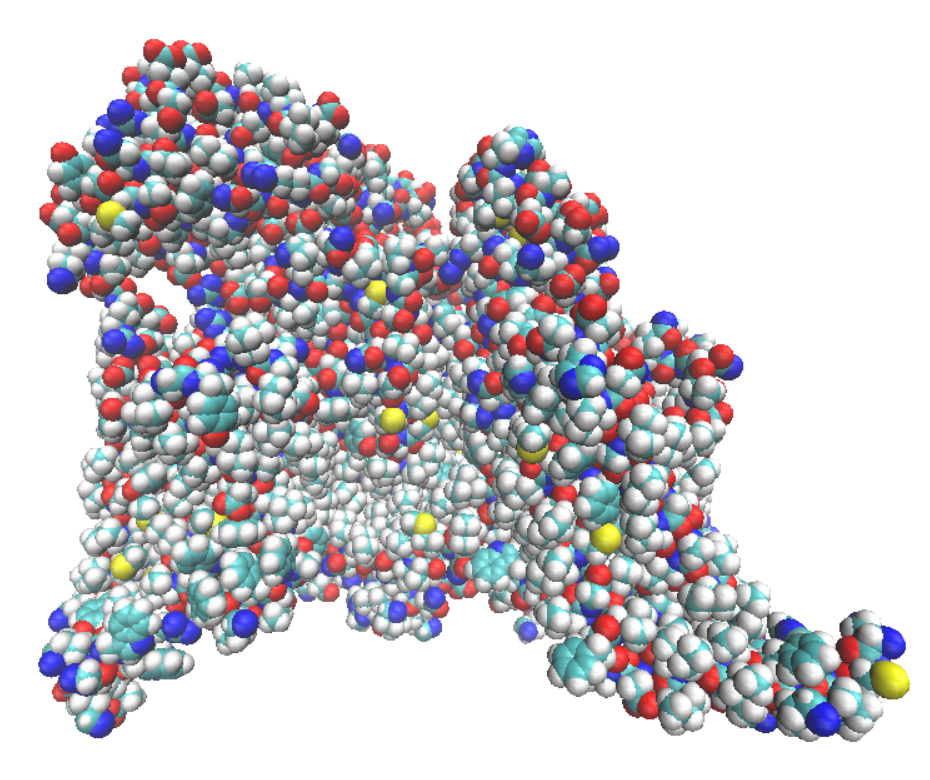}}
  \subfigure[]{
  \label{Gaussian:Surface}
  \includegraphics[scale = 0.23]{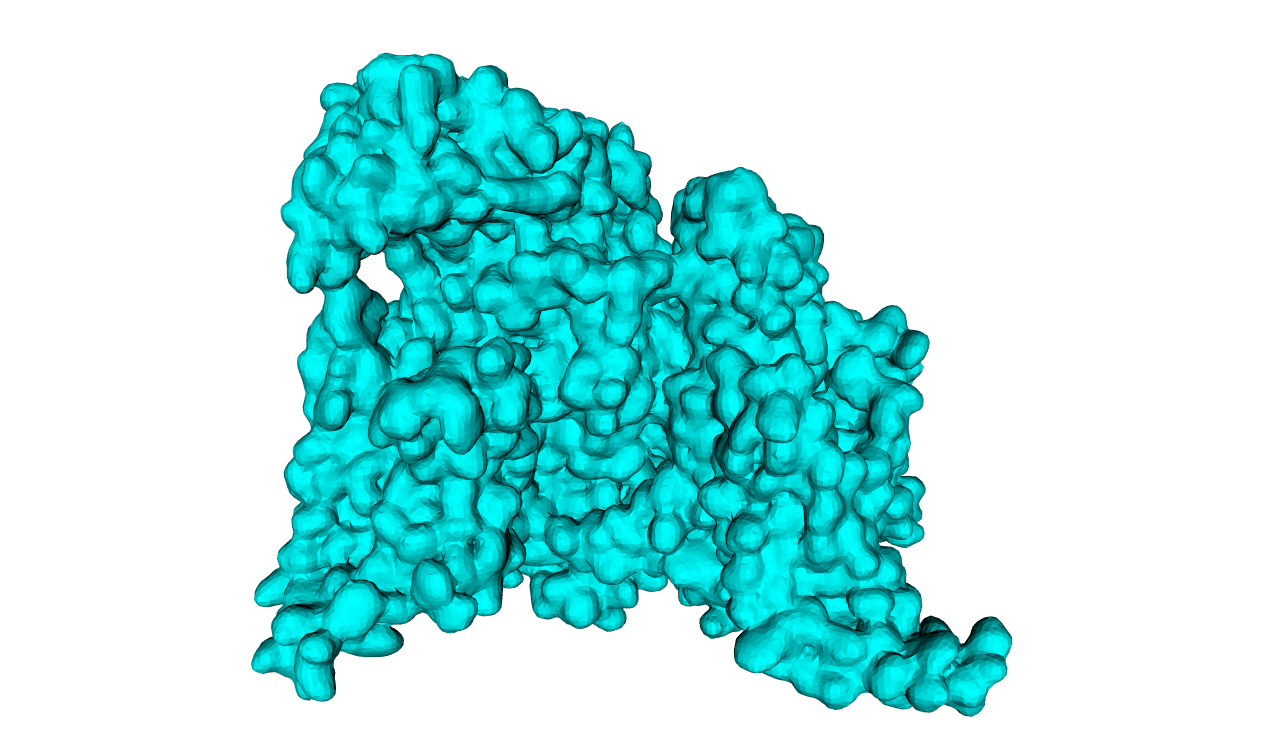}}

  \caption{An example of Gaussian molecular surface(Nav.1.4\cite{Paneaau2486}, PDB code is 6AGF) via VCMM\cite{Bai2014VCMM}. (a) shows the VDW surface, and (b) shows the Gaussian molecular surface generated by TMSmesh \cite{Chen2012Triangulated,Chen2011TMSmesh,Liu2018EFFICIENT} with parameter $d$ and $c$ is 0.9 and 1.0, respectively. All coordinates and corresponding radii are drawn from the PQR file that is transformed from the PDB file, 6AGF, using the PDB2PQR tool \cite{Dolinsky2004PDB2PQR}.}
  \label{GaussianSurface}
\end{figure}

Various definitions of molecular surface exit, including the van der Waals (VDW) surface, solvent accessible surface (SAS) \cite{lee1971the}, solvent excluded surface (SES) \cite{RICHARDS1977AREAS}, molecular skin surface \cite{Edelsbrunner1999Deformable}, the minimal molecular surface \cite{Bates2008Minimal}, and Gaussian surface, etc. The VDW surface is defined as the surface of the union of the spherical atomic surfaces with the VDW radius of each atom in the molecule. The SAS and SES are represented by the trajectory of the center and the interboundary of a rolling probe on the VDW surface, respectively. The molecular skin surface is the envelope of an infinite family of spheres derived from atoms by convex combination and shrinking. The minimal molecular surface is defined as a result of the surface free energy minimization. Different from these definitions, the Gaussian surface is defined as a level set of the summation of the Gaussian kernel functions as follows:
\begin{equation}
\label{GaussLevelSet}
\left\{ \mathbf{x}  \in \mathbb{R}^3,\phi \left( { \mathbf{x}} \right) = c \right\},
\end{equation}
where
\begin{equation}
\label{GaussKernel}
\phi ( \mathbf{x}) = \sum \limits _{i = 1}^N e^{-d(\| \mathbf{x} -  \mathbf{x}_i \|^2 - r_i^2)},
\end{equation}
the parameter $d$ is positive and controls the decay rate of the kernel functions, $\mathbf{x_i}$ and $r_i$ are the location and radius of atom $i$, $c$ is the isovalue, and it controls the volume enclosed by the Gaussian kernel. Comparing with other definitions, the advantages of Gaussian surface present as following: the Gaussian surface is more smooth, and it provides an analytical representation of the electron density of a molecule \cite{DUNCAN1993SHAPE}. The VDW surface, SAS and SES can be approximated well by the Gaussian surface with proper parameter selection \cite{DUNCAN1993SHAPE,Liu2015Parameterization}. The Gaussian surface has been widely used in many problems, for instance, docking problems \cite{mcgann2003gaussian}, molecular shape comparisons \cite{grant1996a}, calculating SAS areas \cite{weiser1999optimization} and the generalized Born models \cite{yu2006what}. Fig. \ref{GaussianSurface} shows an example of a Gaussian surface. This molecule is the structure of the human voltage-gated sodium channel $Na_v 1.4$ in complex with $\beta_1$ \cite{Paneaau2486}. Fig. \ref{Gaussian:PQR} shows all the atoms in the molecule, and Fig. \ref{Gaussian:Surface} shows the corresponding Gaussian surface.

For traditional Gaussian surface, the level set function is constructed by a summation of Gaussian kernel functions, whose number depends on the total number of atoms in the molecular. For large biomolecules, the number of kernels in their definition of Gaussian molecular surface may achieve millions. It leads to a significant challenge for their analysis and understanding. Zhang et al. \cite{Liao2015Atom} proposed an atom simplification method for the biomolecular structure based on Gaussian molecular surface. This method contains two main steps. The first step eliminates the low-contributing atoms. The second step optimizes the center location, the radius and the decay rate of the remaining atoms.

In the field of computer aided geometric design, the Gaussian surface is a typical implicit surface representation method. In the past several decades, there are wealthy works focusing on the implicit surface reconstruction problem, and various approaches have been proposed. J.C. Carr \cite{Carr2001Reconstruction} employed a method to reconstruct an implicit surface with RBFs based on a greedy algorithm to append centers with large residuals to decrease the number of basis functions. However, the solution of this method is not sparse enough. M. Samozino \cite{Samozino2006RVC} put forward a strategy to put the RBFs’ centers on the Voronoi vertices. This method, firstly chooses a user-specified number of centers by filtering and clustering from a subset of Voronoi vertices, secondly gets the reconstructed surface by solving a least-square problem. Unfortunately, it causes larger approximation error on the surface while approximating the surface and center points equally. In 2016, Chen \cite{2016LiSparse} et al. proposed a model of sparse RBF surface representations. They constructed the implicit surface based on sparse optimization with general RBF. And the initial Gaussian RBF is on the medial axis of the input object. They have solved the RBF surface by sparse optimization technique. Sparse optimization has become a very popular technique in many active fields, such as signal processing and computer vision, etc \cite{Elad2010}. This technique has been applied in linear regression \cite{Atkinson1992Subset}, deconvolution \cite{TAYLOR1979DECONVOLUTION}, signal modeling \cite{RISSANEN1978MODELING}, preconditioning \cite{Grote1997Parallel}, machine learning \cite{Girosi1998equivalence}, denoising \cite{Chen1998Atomic}, and regularization \cite{Daubechies2004sparsity}. In recent years, sparse optimization also has been applied in geometric modeling and graphics problems (refer to a review \cite{Xu2015}).

In this paper, based on the frame of sparse optimization, we propose a model for reducing the number of kernels in the definition of Gaussian surface while preserving the shape of the molecular surface. We emphasize several differences between our method and previous sparse optimization methods with surface representation: 1) Compared with other works, our focus is mainly on reducing the number of kernels in Gaussian molecular surface; 2) The objective function of our model is a complicated nonlinear function with respect to the locations, sizes, shapes and orientations of RBFs; 3) Different initializations and algorithms are proposed for solving the corresponding sparse optimization problem in our model.

The outline of this paper is as follows. Section \ref{section2} reviews some preliminary knowledge about radial basis functions and sparse optimization. Section \ref{section3} presents our model together with an algorithm for representing the Gaussian molecular surface sparsely. The experimental results and comparisons are demonstrated in section \ref{section4}. We conclude the paper in section \ref{section5}.

\section{Preliminaries}\label{section2}
In this section, we present some preliminary knowledge for radial basis functions, ellipsoid Gaussian RBF and sparse optimization.
\subsection{Radial basis functions}
Consider an implicit function, $\phi(\mathbf{x})=const$, where $\phi(\mathbf{x})$ is expressed by a linear combination of some basis functions as follows:
\begin{equation}
\label{RFBdefined}
\phi(\mathbf{x}):= \sum \limits _{i = 1}^n  c_i \xi _i(\mathbf{x}),
\end{equation}
where $\xi _i (\mathbf{x}), i = 1,2,\cdots, n$ are basis functions, and $c_i, i = 1,2,\cdots,n$ are combination coefficients. $\mathbf{x}$ is the variables in the geometric space, for instance, in two dimensional space $\mathbf{x} = (x_1,x_2)^\top$, and in three dimensional space $\mathbf{x} = (x_1,x_2,x_3)^\top$.

In this paper, the RBF is rewritten as $\xi _i (\mathbf{x}) = \xi (\|\mathbf{x} - \mathbf{c}_i\|)$, where $\xi (\sigma)$ is a nonnegative function defined on $[0,+ \infty)$ which generally has the following properties \cite{Carr2001Reconstruction}:

\begin{enumerate}
  \item $\xi (0) = 1$;
  \item $lim _{\sigma \rightarrow + \infty} \xi (\sigma) = 0$;
\end{enumerate}
$\textbf{c}_i$ is the center of the basis function $\xi _i (\textbf{X})$,

A typical choice of radial basis function is Gaussian function
\begin{equation}
\xi (\sigma) = e^{-d\sigma ^2}
\end{equation}
where $d \in \mathbb{R}$ is a decay value.
Besides the Gaussian function, there are other radial basis functions including inverse multiquadric radial basis functions, for instance, $\xi (\sigma) = \frac{1}{\sqrt{1 + \sigma ^2}}$ for $\sigma \in \mathbb{R}$.

\subsection{Ellipsoid Gaussian RBF}
\label{MGRBF}
In general, Gaussian RBF is defined as:
\begin{equation}
\label{GenGaussRBF}
\xi _g (\mathbf{y}) = e^{-\|\mathbf{y} - \mathbf{x}\|^2} ,
\end{equation}
where $\mathbf{x}$ is the Gaussian RBF center. Let $\xi _g (\mathbf{y}) = const$, and we can obtain a general sphere equation as follows:
\begin{equation}
\label{sphereEq}
\begin{aligned}
\| \mathbf{y} - \mathbf{x}\|^2 = const.\\
\end{aligned}
\end{equation}
We only consider the case of $\mathbf{y} = (y_1,y_2,y_3)^\top \in \mathbb{R}^3$, and $\mathbf{x} = (x_1,x_2,x_3)^\top\in \mathbb{R}^3$, then Eq. \ref{sphereEq} is rewritten as:
\begin{equation}
\label{SphereEq}
\left(y_{1}-x_{1}\right)^{2}+\left(y_{2}-x_{2}\right)^{2}+\left(y_{3}-x_{3}\right)^{2}=const,
\end{equation}
Eq. \ref{SphereEq} can be written as a quadratic form and is equivalent to Eq. \ref{quadratic}:
\begin{equation}
\label{quadratic}
(\mathbf{y}-\mathbf{x})^\top \left[
\begin{array}{lll}
{1} & {0} & {0} \\
{0} & {1} & {0} \\
{0} & {0} & {1}
\end{array}\right](\mathbf{y}-\mathbf{x})=const.
\end{equation}
Thus, we regard the level set of the general Gaussian RBF function as the sphere equation in the 3D Euclid space. However, for complex geometric shapes, the general Gaussian RBF has certain limitations to approximation \cite{2016LiSparse}, which also leads to requirement of a large number of basis functions while approaching the input surface. Thus, we propose a new ellipsoid Gaussian RBF.

Firstly, we can modify the standard sphere equation to the standard ellipsoid equation. We just introduce three coefficient for the x,y,z direction, the form is as follows:
\begin{equation}
\label{scalablepll}
d_1 \left(y_{1}-x_{1}\right)^{2}+d_2 \left(y_{2}-x_{2}\right)^{2}+d_3 \left (y_{3}-x_{3}\right)^{2}=const.
\end{equation}
Obviously, Eq. \ref{SphereEq} is the special case of Eq. \ref{scalablepll} when $d_1 = d_2 = d_3$. Likewise, Eq. \ref{scalablepll} can be rewritten as a quadratic form: 
\begin{equation}
\label{epllquad}
(\mathbf{y}-\mathbf{x})^\top \left[
\begin{array}{lll}
{d_{1}} & {0} & {0} \\
{0} & {d_{2}} & {0} \\
{0} & {0} & {d_{3}}
\end{array}\right](\mathbf{y}-\mathbf{x})=const,
\end{equation}
and we can obtain the standard ellipsoid Gaussian RBF:

\begin{equation}
\label{scalableRBF}
\xi _t (\mathbf{y}) = e^{-(\mathbf{y} - \mathbf{x})^\top \mathbf{D}(\mathbf{y} - \mathbf{x})},
\end{equation}
where $\mathbf{D} = diag(d_{1},d_{2},d_{3})$, $d_1, d_2, d_3 \in \mathbb{R}$.

In the next stage, we construct the general ellipsoid Gaussian RBF. The method introduces rotation mechanism for the ellipsoid feature. The form is as follows:
\begin{equation}
(\mathbf{y}-\mathbf{x})^\top \mathbf{R}^\top \left[
\begin{array}{lll}
{d_{1}} & {0} & {0} \\
{0} & {d_{2}} & {0} \\
{0} & {0} & {d_{3}}
\end{array}\right]\mathbf{R}(\mathbf{y}-\mathbf{x})=const,
\end{equation}
where $\mathbf{R}$ is the total rotation matrix, and it is equal to the product of rotation matrices from three directions.
\begin{equation}
\label{RotataionMatirx}
\mathbf{R}(\alpha, \beta, \gamma) = \mathbf{R_z}(\gamma ) \cdot \mathbf{R_y}(\beta )\cdot \mathbf{R_x}(\alpha ),
\end{equation}
and $\mathbf{R_x}(\alpha )$ is a rotation matrix of x direction:
\begin{equation}
\mathbf{R_{x}}\left(\alpha \right)=\left[ \begin{array}{ccc}{1} & {0} & {0} \\ {0} & {\cos \alpha} & {-\sin \alpha} \\ {0} & {\sin \alpha} & {\cos \alpha}\end{array}\right],
\end{equation}
$\mathbf{R_y(\beta )}$ is a rotation matrix of y direction:
\begin{equation}
\mathbf{R_{y}}\left(\beta\right)=\left[ \begin{array}{ccc}{\cos \beta} & {0} & {-\sin \beta} \\ {0} & {1} & {0} \\ {\sin \beta} & {0} & {\cos \beta}\end{array}\right],
\end{equation}
$\mathbf{R_z(\gamma )}$ is a rotation matrix of z direction:
\begin{equation}
\mathbf{R_{z}}\left(\gamma\right)=\left[ \begin{array}{ccc}{\cos \gamma} & {-\sin \gamma} & {0} \\ {\sin \gamma} & {\cos \gamma} & {0} \\ {0} & {0} & {1}\end{array}\right],
\end{equation}
so that $\mathbf{R(\alpha , \beta , \gamma )}$ is equal to :
\begin{equation}
\left[
\begin{array}{ccc}
{\cos \beta \cos \gamma} & {-\cos \alpha \sin \gamma -  \sin \alpha \sin \beta \cos \gamma} & {\sin \alpha \sin \gamma - \cos \alpha \cos \gamma \sin \beta} \\
{\cos \beta \sin \gamma} & {\cos \alpha \cos \gamma - \sin \alpha \sin \beta \sin \gamma} & {-  \sin \alpha \cos \gamma - \cos \alpha \sin \beta \sin \gamma} \\
{\sin \beta } & {\cos \beta \sin \alpha} & {\cos \alpha \cos \beta}
\end{array}
\right].
\end{equation}

Then we can construct an ellipsoid Gaussian RBF, and the equation is as follows:
\begin{equation}
\label{RAndTRBF}
\widetilde{\xi} (\mathbf{y}) =  e^{-\|\mathbf{D}^{\frac{1}{2}}\mathbf{R}(\mathbf{y} - \mathbf{x})\|_2^2}.
\end{equation}

The summation of ellipsoid Gaussian RBF can be rewritten as:
\begin{equation}
\label{ImplicitModel}
\tilde{\phi} (\mathbf{y}) = \sum \limits _{i = 1}^N c_i \widetilde{\xi}_i(\mathbf{y}) = \sum \limits _{i = 1}^N c_ie^{-\|\mathbf{D}_i^{\frac{1}{2}}\mathbf{R}_i(\mathbf{y} - \mathbf{x}_i)\|_2^2}.
\end{equation}


\subsection{Sparse optimization}
 $\mathbf{x} = (x_1,x_2,\cdots,x_n)^\top$ is an n-dimensional vector. The $L_{p}$ norm of $\mathbf{x}$ is defined as $\|\mathbf{x}\|_{p} = (\sum _{i = 1}^n |x_i|^p)^{1/p} (0 < p < + \infty)$. The $L_0$  norm of $\mathbf{x}$ is defined as $\|\mathbf{x}\|_{0} :=\#\left\{x_{i} \| x_{i} \neq 0\right\}$, i.e., the number of nonzero elements in $\mathbf{x}$. If $\|\mathbf{x}\|_0 \ll n$, the vector $\mathbf{x}$ is called sparse.

Consider solving a linear system of equations
\begin{equation}\label{overdeterminedSystem}
\mathbf{A}\mathbf{x} = \mathbf{b},
\end{equation}
where the rank of $\mathbf{A}$ is much less than the number of columns of $\mathbf{A}$.

Assuming there are infinite number of solutions in the system defined by Eq. \ref{overdeterminedSystem}, our goal is to find the sparsest solution among all the feasible solutions. This is equivalent to solving the following optimization problem

\begin{equation}
\label{L0}
\begin{aligned}
& \min _{\mathbf{x}} \ \ \| \mathbf{x} \|_0 \\
& s.t. \mathbf{A}\mathbf{x} = \mathbf{b}.
\end{aligned}
\end{equation}

The $L_0$ optimization is a NP-hard problem. To solve Eq. \ref{L0} directly, one must sift through all possible distributions of the nonzero components in $\mathbf{x}$. This method is intractable because the search space is exponentially large \cite{Rao1999BBS,Natarajan1995Sparse}. To surmount this obstacle, one might replace the $L_0$ quasinorm with the $L_1$ norm to obtain a convex optimization problem

\begin{equation}
\label{L1}
\begin{aligned}
& \min _{\mathbf{x}} \ \ \| \mathbf{x} \|_1 \\
& s.t. \mathbf{A}\mathbf{x} = \mathbf{b}.
\end{aligned}
\end{equation}

Intuitively, the $L_1$ norm is the convex function closest to the $L_0$ quasi-norm, so this substitution is referred to convex relaxation. It is shown that the solutions of the two problems (Eq. \ref{L0}) and (Eq. \ref{L1}) are equivalent under certain conditions\cite{Troop2006Relax}. Problem defined in Eq. \ref{L1} is convex and can be efficiently solved.

The problem (Eq. \ref{L1}) can be rewirtten as an unconstrained optimization problem

\begin{equation}
\label{NewL1}
\min _{\mathbf{x}} \ \ \alpha \| \mathbf{x} \|_1 + \beta \|\mathbf{A}\mathbf{x} - \mathbf{b}\|_2^2.
\end{equation}
where $\alpha > 0$ and $\beta > 0$ are parameters which balances the two targets: sparsity and accuracy of solutions.

Also, Eq. \ref{NewL1} is called the least absolute shrinkage and selection operator(LASSO) and can be solved by many well established algorithms, such as ADMM, CG, LBFGS.

\section{Model and algorithm}\label{section3}
In this section, base on sparse optimization frame, we present a new optimization model for finding a function $\tilde{\phi}$ having fewest possible basis with shape of ellipsoid to approximate the summation of Gaussian kernel functions $\phi$, in original Gaussian molecular surface. Then, we propose the numerical optimization algorithm to solve this model.
\subsection{Modeling with ellipsoid Gaussian RBF}

According to the definition of Gaussian surface, we first give the model of representing the Gaussian molecular surface sparsely as follows:
\begin{equation}\label{model1}
\begin{array}{c}
\mathop {\min }\limits_{\mathbf{X}} {w_s} \cdot E_s(\mathbf{X}) + {w_l} \cdot E_{l1}\left(\mathbf{c}, \mathbf{d}_{p}\right)\\
s.t.\left\{ \begin{array}{ll}
c_i  \ge 0 \qquad &i = 1,2,\cdots,N, \\
d_{pi}  \ge 0 \qquad &p = 1,2,3.
\end{array} \right.
\end{array}
\end{equation}
$E_s(\mathbf{X})$ is the error between $\tilde{\phi}$ and $\phi$  at constrained points $\mathbf{y}_k,k=1,2, \cdots, M$. $\tilde{\phi}$ is a summation of ellipsoid Gaussian RBFs. $\phi(\mathbf{x})$ is the implicit function in the definition of Gaussian molecular surface (Eq. \ref{GaussKernel}). It is to be approximated by $\tilde{\phi}(\mathbf{x})$.
\begin{equation}
E_s(\mathbf{X}) = \sum \limits _{k=1}^{M}\left[\tilde{\phi}(\mathbf{y_k};\mathbf{X}) - \phi(\mathbf{y_k})\right]^{2}  = \sum \limits _{k=1}^{M}\left[ \sum \limits _{i = 1}^N c_i  e^{-\|\mathbf{D}_i^{\frac{1}{2}}\mathbf{R}_i(\alpha_i,\beta_i,\gamma_i)(\mathbf{y}_k - \mathbf{x}_i)\|_2^2} - \phi(\mathbf{y_k})\right]^{2},
\end{equation}
where $\mathbf{y}_k = (y_{k1},y_{k2},y_{k3})^\top \in \mathbb{R}^3$  is the $k$th constrained point. $\mathbf{x}_i = (x_{i1},x_{i2},x_{i3})^\top\in \mathbb{R}^3$ is the center of the $i$th ellipsoid Gaussian RBF, $\mathbf{D}_i = diag(d_{i1},d_{i2},d_{i3})$ define the lengths of ellipsoid along three main axis, $\mathbf{R}_i(\alpha_i,\beta_i,\gamma_i)$ is a rotation matrix, the $\alpha_i,\beta_i,\gamma_i$ are rotation angles of the $i$th ellipsoid Gaussian RBF, $i = 1,2,\cdots, N$, $k = 1,2,\cdots,M$. $N$ is the number of the centers. $M$ is the number of the constrained points. The $\mathbf{X}$ is the optimization variable, the formula is
\begin{equation}\label{X}
\mathbf{X} = \left[\mathbf{c}, \mathbf{d}_{p}, \mathbf{x}, \bm{\alpha}, \bm{\beta}, \bm{\gamma} \right]^\top, i = 1,2,\cdots,N, p = 1,2,3,
\end{equation}
where $\mathbf{c}=\left[c_{1,} c_{2}, \cdots, c_{N}\right]$, $\mathbf{d}_{p}=\left[d_{1p}, d_{2p}, \cdots, d_{N p}\right]$, $\mathbf{x}=\left[\mathbf{x}_{1}, \mathbf{x}_{2}, \cdots, \mathbf{x}_{N}\right]$, $\bm{\alpha}=\left[\alpha_{1,} \alpha_{2}, \cdots, \alpha_{N}\right]$, $\bm{\beta}=\left[\beta_{1,} \beta_{2}, \cdots, \beta_{N}\right]$, $\bm{\gamma}=\left[\gamma_{1,} \gamma_{2}, \cdots, \gamma_{N}\right]$.


The second term of objective function is the sparse $L_1$ norm of $\mathbf{d}_p$ and $\mathbf{c}$.
\begin{equation}\label{Sparseterm}
E_{l1}(\mathbf{c},\mathbf{d}_p) = \|\mathbf{c}\|_1 + \sum \limits _{p = 1}^3 \|\mathbf{d}_p\|_1 = \sum \limits _{i = 1}^N |c_i| + \sum \limits _{i = 1}^N |d_{i1}| + \sum \limits _{i = 1}^N |d_{i2}| + \sum \limits _{i = 1}^N |d_{i3}|.
\end{equation}

In Chen's work\cite{2016LiSparse}, they use general Gaussian RBF to represent an implicit surface with sparse optimization model. Their fitting function is a linear combination of general Gaussian RBF. Consequently, their sparse term is just the $L_1$ norm of coefficient of basis functions. Based on our model, especially for the ellipsoid Gaussian RBF, in order to find fewer basis functions, we give corresponding modified sparse term in Eq. \ref{Sparseterm}. Our model is equivalent to a LASSO regression, in which the locations, decay parameters and directions of ellipsoid Gaussian RBFs are changeable to achieve low errors at constrained points and high sparsity of $\mathbf{c}$ and $\mathbf{d}_p$.

The $w_s > 0$ and $w_l > 0$ are parameters which balances the two targets: accuracy of solutions $E_s$ and sparsity $E_{l1}$. And the constrained conditions are explained as follows.
\begin{itemize}
  \item $c_i>0$ indicates that the corresponding ellipsoid Gaussian RBF is nonnegative which means each RBF in $\tilde{\phi}$ can be seen as a new real physical atom with ellipsoid shape.
  \item $d_{pi} \geq 0$ implies  the basis function is zero at infinity, which is consistent with the fitted function $\phi$.
\end{itemize}

In order to transform the Eq. \ref{model1} to an unconstrained optimization problem, we do the following substitution,
\begin{equation}
\left\{ \begin{array}{ll}
c_i  = \tilde{c_i}^2 \qquad &i = 1,2,\cdots,N \\
d_{pi}  = \tilde{d}_{pi}^2 \qquad &p = 1,2,3,\\
\end{array} \right.
\end{equation}
And corresponding $\mathbf{\tilde{D}}_i = diag(\tilde{d}_{i1}^2,\tilde{d}_{i2}^2,\tilde{d}_{i3}^2)$. For simplicity, we still use $c_i,d_{pi},\mathbf{D}_i$ to denote $\tilde{c}_i,\tilde{d}_{pi},\mathbf{\tilde{D}}_i$.
Thus, the equivalent unconstrained model is:
\begin{equation}\label{Model2}
\mathop {\min }\limits_{\mathbf{X}} {w_s} \cdot \sum_{k=1}^{M}\left[\sum_{i=1}^{N} c_i^2 \cdot e^{-\| \mathbf{D}_i\mathbf{R}_i(\mathbf{y}_{k} - \mathbf{x}_i)\|^2_2} - \phi\left(\mathbf{y}_{k}\right)\right]^{2} + w_l \cdot \left\{\|\mathbf{c}\|_1 +\sum \limits _{p = 1}^3 \|\mathbf{d}_p\|_1 \right\}.\\
\end{equation}
\subsection{Optimization algorithm}
\subsubsection{Overview}
\begin{figure}[H]
  \centering
  \includegraphics[scale = 0.26]{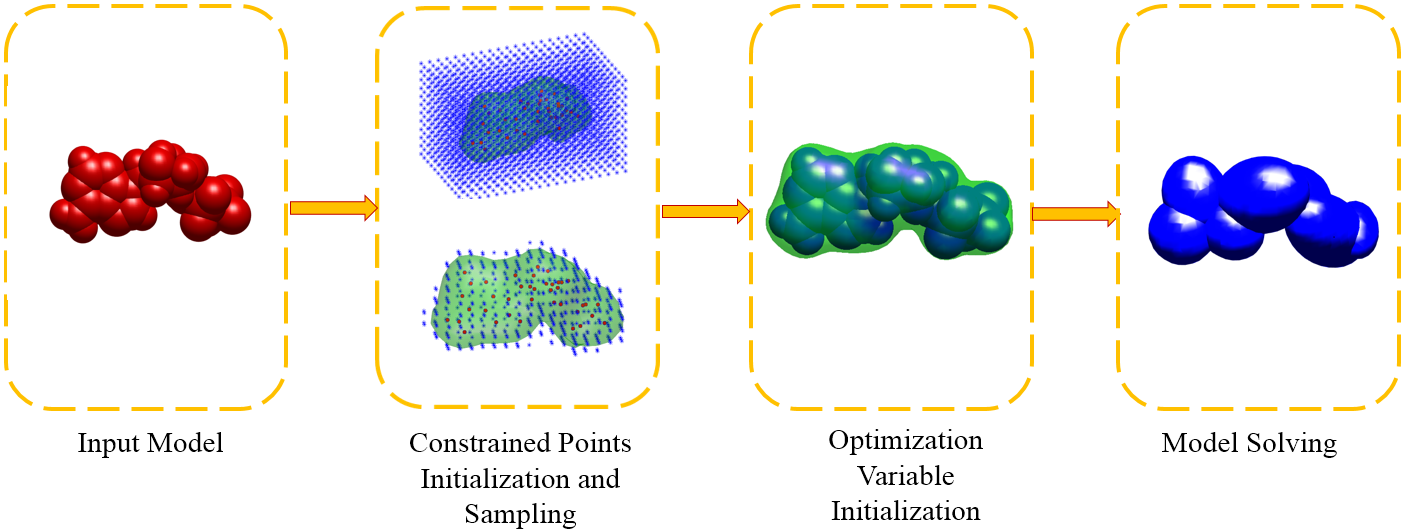}\\
  \caption{The process of our algorithm and results in each step.}\label{OverView}
\end{figure}

In this section, we present algorithms for solving the minimization problem defined in Eq. \ref{Model2}. The process of solving Eq. \ref{Model2} consists of three parts. In the first part, the constrained points $\mathbf{y}_k, k=1,...,M$ need to be set.  They are selected from vertices of a three dimensional orthogonal grid bounding all atoms and they are closed to the original Gaussian surface, $\phi(\mathbf{x})=c$; In the second part, we set the proper initial values of $\mathbf{X}$ in Eq. \ref{Model2}; In the third part, we solve the problem defined in Eq. \ref{Model2} using an algorithm of minimizing the sparsity and error terms in Eq. \ref{Model2} alternatively. Fig. \ref{OverView} demonstrates the process of our algorithm. The result shows that, using our method, the original Gaussian surface is approximated well by a summation of much fewer ellipsoid Gaussian RBFs.

\subsubsection{Constrained points initialization and sampling}
\label{Yinitialization}
\begin{figure}[H]
  \centering
  \subfigure[]{
  \label{point:surface}
  \includegraphics[scale = 0.38]{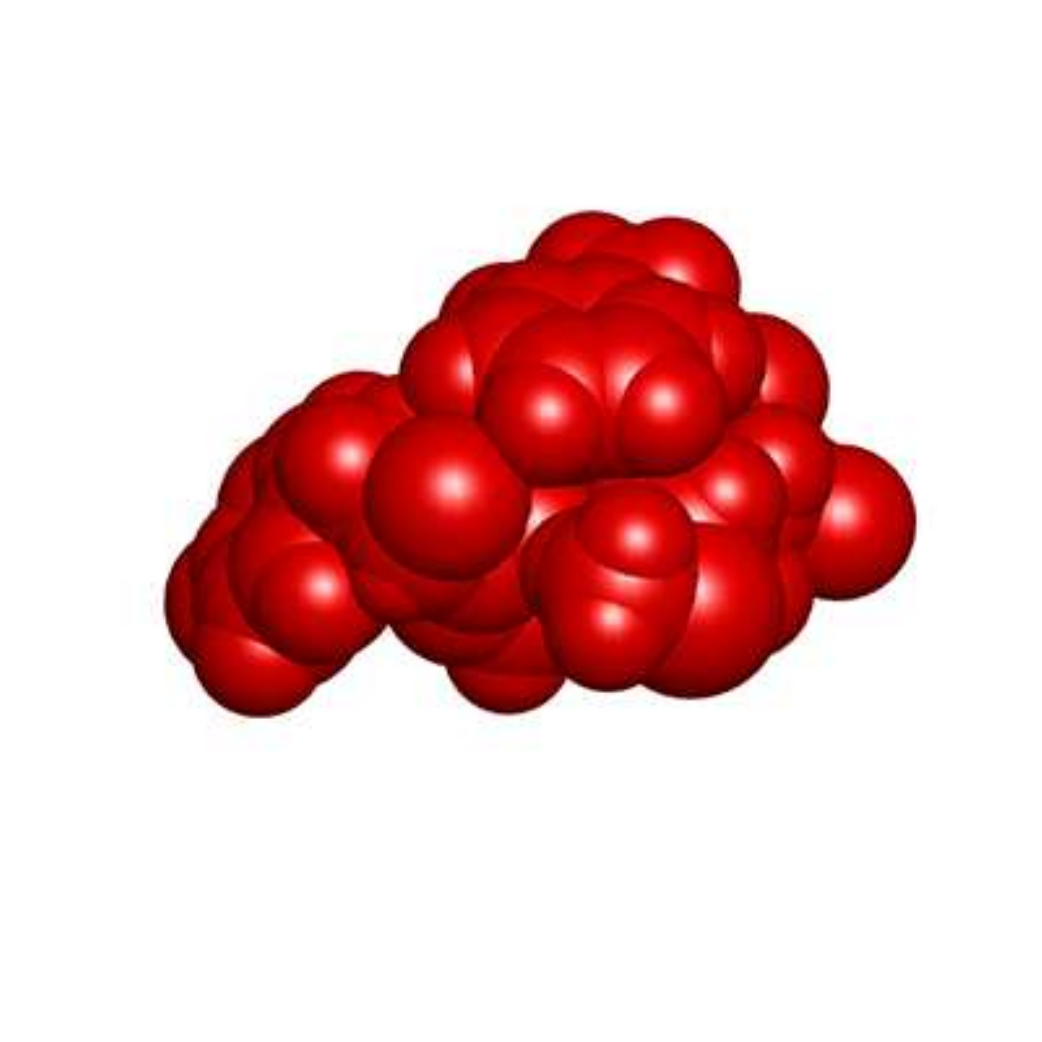}}
  \subfigure[]{
  \label{point:box}
  \includegraphics[scale = 0.38]{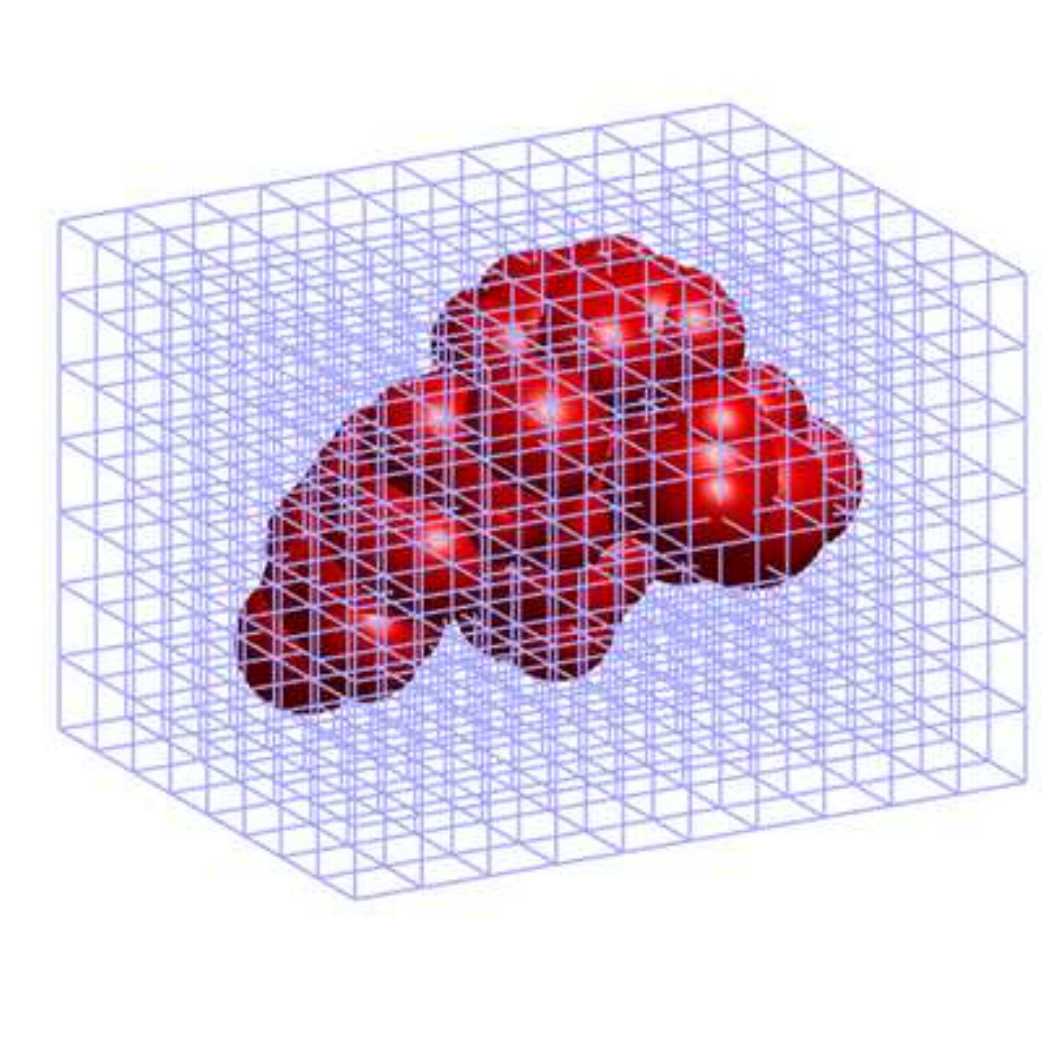}}
  \subfigure[]{
  \label{point:constraint}
  \includegraphics[scale = 0.38]{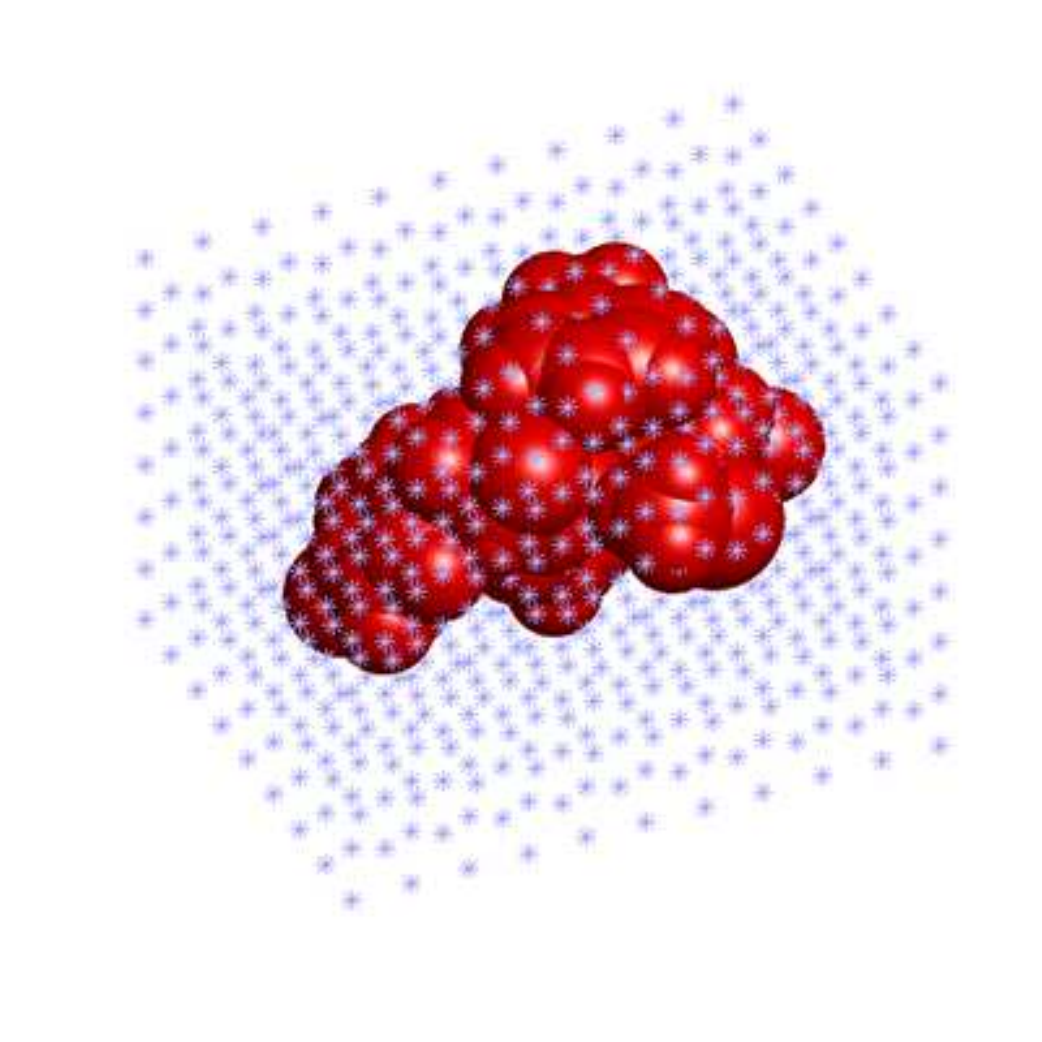}}
  \caption{Constrained points initialization. (a) shows a real molecule (PDBID: 2LWC). (b) shows the bounding box of the molecular. (c) shows constrained points.}
  \label{point}
\end{figure}

For simplicity, we put the molecule (Fig. \ref{point:surface}) on the bounding box $\Omega$ (Fig. \ref{point:box}) in 3D space. The range of bounding box is $[a,b] \times [c,d] \times [e,f]$, where $a,b,c,d,e,f \in \mathbb{R}$.
The bounding box $\Omega$ is discretized into a set of uniform grid points as shown in Fig. \ref{point:constraint}:
\begin{equation}
\left\{\mathbf{P}_{ijk}\right\} = \left\{(x_{i}, y_{j}, z_{k})\right\}, \quad  i=1, 2,  \cdots, N_{x}, j=1, 2, \cdots, N_{y}, k=1, 2, \cdots, N_{z},
\end{equation}
where 
\begin{equation}
\begin{aligned}
x_i &= a + i \times \frac{(b - a)}{N_x}, y_j &= c + j \times \frac{(d - c)}{N_y}, z_k &= e + k \times \frac{(f - e)}{N_z},
\end{aligned}
\end{equation}

In order to reduce the scale of optimization problem and make the problem solvable in usual machine memory and time, it is necessary to decrease the number of constrained points. A set of constrained points $\left\{\mathbf{y}_{k}\right\}_{k=1}^{M}$ is selected from the set of grid points $\left\{\mathbf{P}_{ijk}\right\}$, and $\left\{\mathbf{y}_{k}\right\}_{k=1}^{M}$ are close to the original Gaussian surface $\phi(\mathbf{x})=c$. In our practice, the constrained points $\left\{\mathbf{y}_{k}\right\}_{k=1}^{M}$ stratifying $\|\phi(\mathbf{y}_k)-c\| \leq 1$ are chosen.


\subsubsection{Optimization variable initialization}
\label{Vinitialization}
In this section, we give the initialization of optimization variable $\mathbf{X}$ defined in Eq. \ref{X}, consisting of $\mathbf{c}$, $\mathbf{d}_p$, $\mathbf{x}$, $\bm{\alpha}$, $\bm{\beta}$, $\bm{\gamma}$, where $\mathbf{c}$ are the coefficients of ellipsoid RBFs, $\mathbf{d}_p$ are decay rates of ellipsoid Gaussian RBFs, $\mathbf{x}$ are center coordinates of ellipsoid Gaussian RBFs, $\bm{\alpha},\bm{\beta},\bm{\gamma}$ are rotate angles of ellipsoid Gaussian RBFs. Since $\phi(\mathbf{x})$ is given by a summation of  Gaussian RBFs located at given atoms and general Gaussian RBF is a degradation case of ellipsoid Gaussian RBF, $\tilde{\phi}$ can be initialized as the same as $\phi$. \\
\textbf{The decay rates $\mathbf{d}_p$ and rotate angles $\bm{\alpha},\bm{\beta},\bm{\gamma}$  of ellipsoid Gaussian RBFs}

In the definition of Gaussian molecular surface (Eq. \ref{GaussKernel}) the value of d is usually set to be 0.5. So the decay rates $\mathbf{d}_p$ in $\tilde{\phi}$ is also set to be a constant vector as follows
\begin{equation}
\mathbf{d}_{p} = [\underbrace{0.5,0.5,\cdots,0.5}_{N}]^{\top}.
\end{equation}
The initial angles in $\tilde{\phi}$ are set to be zeros,
\begin{equation}\label{angle}
\bm{\alpha} = \mathbf{0},\bm{\beta} = \mathbf{0},\bm{\gamma} = \mathbf{0}.
\end{equation}
\textbf{The center coordinates $\mathbf{x}$ of ellipsoid Gaussian RBFs}

The traditional methods set the centers of the RBFs on the input surface and/or offset points. Chen \cite{2016LiSparse} et al present a way of putting the centers of the RBFs on the medial axis of the input object. Different from their work, we address our model and method on molecular Gaussian surface and the initial values of centers of ellipsoid Gaussian RBFs are directly given by the centers of atoms as follows.
\begin{equation}\label{center}
\mathbf{x}_{i} = \left[x_{atom}^{(i)},y_{atom}^{(i)},z_{atom}^{(i)}\right]^{\top}, i = 1,2, \cdots,N,
\end{equation}
where $x_{atom}^{(i)},y_{atom}^{(i)},z_{atom}^{(i)}$ are coordinates of $i$th atom.\\
%
\textbf{The coefficients $\mathbf{c}$ of ellipsoid Gaussian RBFs}

Once $\mathbf{d}_p,\mathbf{x}_i,\bm{\alpha},\bm{\beta},\bm{\gamma}$ have been chosen and atom radii $\mathbf{r} = \left[r_1,r_2,\cdots,r_N\right]^{\top}$  is given, to initialize $\tilde{\phi}$ as the same as $\phi$, we set the coefficients $\mathbf{c}$ of  ellipsoid Gaussian RBF as following.
\begin{equation}\label{c_init}
\mathbf{c} = \left[\sqrt{e^{dr_1^2}},\sqrt{e^{dr_2^2}},\cdots,\sqrt{e^{dr_N^2}}\right]^{\top}.
\end{equation}

\subsubsection{Sparse optimization}
After initialization of problem Eq. \ref{Model2}, we hope to find a sparse coefficients vector $\mathbf{c} = [c_1, c_2, \cdots, c_N]^{\top}$ and a sparse decay rates vector $\mathbf{d} = [\mathbf{d}_1, \mathbf{d}_2, \mathbf{d}_3]^{\top}$  of $\tilde{\phi}$($\mathbf{X})$ as defined in Eq. \ref{ImplicitModel} such that $\tilde{\phi}$ is a good fit of $\phi$ at the constrained points $\left\{\mathbf{y}_{k}\right\}_{k=1}^{M}$. Algorithm \ref{frame} represents the main modules of our sparse optimization method, which is described below.
\begin{algorithm}
\setstretch{1.05}
\caption{Sparse optimization}
\label{frame}
\begin{algorithmic}[1]
\STATE{Input: PQR file containing coordinates of centers and radii of atoms}
\STATE{Output: The list of parameters of ellipsoid RBFs, i.e. solution of $\mathbf{X}$ in problem Eq. \ref{Model2}}
\STATE{\textbf{Step 1}. initialize $\mathbf{X}$ as shown in Section \ref{Vinitialization}  }
\STATE{\textbf{Step 2}. select  constrained points $\left\{\mathbf{y}_{k}\right\}_{k=1}^{M}$ as shown in Section \ref{Yinitialization}}
\STATE{\textbf{Step 3}. set the number of maximum iteration $MaxNiter$ and number of sparse optimization iteration $SparseNiter$}
\STATE{\textbf{Step 4}. initialize variable of iteration: $Niter= 0$ and set tolerance: $tol = 1e-3$}
\STATE{\textbf{Step 5}. optimization of $E_s$ and $E_{l1}$ defined in Eq. \ref{Model2} alternatively}
\WHILE{$Niter < MaxNiter$}
\STATE {$Niter = Niter+ 1$}
\STATE {\textbf{Step 5.1}. delete useless Gaussian basis function $|c_i| < tol$ every 20 steps}
\STATE {\textbf{Step 5.2}. calculate $\tilde{\phi}(\mathbf{y}_{k})$ for all constrained points by $\mathbf{X}$}
\STATE \textbf{Step 5.3}. calculate the accuracy term $E_{s}$ and sparse term $E_{l1}$
\STATE \textbf{Step 5.4}. calculate the adaptive coefficients $w_s$ and $w_{l}$
\STATE $\quad\quad\quad\quad\quad\quad$ $w_{s}=\max \left\{\frac{E_{s}}{E_{s}+E_{l 1}}, \varepsilon\right\}$, $w_{l}=\frac{E_{l 1}}{E_{s}+E_{l 1}}$.
\STATE \textbf{Step 5.5}. check the maxium of error between $\tilde{\phi}$ and $\phi$ at constrained points $\mathbf{y}_k$ and correct the coefficients $w_s$ and $w_{l}$
\IF {$\mathop {\max }   _{1 \leq k \leq M} \| \tilde{\phi}(\mathbf{y}_k) - \phi(\mathbf{y}_k) \| > 0.5$}
\STATE $\quad\quad$ $w_s = 1, w_l = 0$
\ENDIF
\STATE \textbf{Step 5.6}. accucacy optimiztion for $E_s$ by set coefficients $w_s$ and $w_{l}$
\IF {$Niter > SparseNiter$}
\STATE $\quad\quad$ $w_s = 1, w_l = 0$
\ENDIF
\STATE \textbf{Step 5.7}. calculate the gradient of object function
\STATE \quad\quad\quad\quad\quad\quad\quad\quad $\nabla f = \nabla \cdot (w_s \cdot E_{s} +  w_{l} \cdot E_{l1})$
\STATE \textbf{Step 5.8}. $\mathbf{X}$ is updated by $X=X-\tau\nabla f$, $\tau$ is computed by inexact line search along the direction of $-\nabla f$
\ENDWHILE
\end{algorithmic}
\end{algorithm}

Step 1 shows initialization of optimization variable $\mathbf{X}$. Step 2 selects constrained points $\left\{\mathbf{y}_{k}\right\}_{k=1}^{M}$. Step 3 and step 4 initialize some variables, i.e. the number of total iterations, the number of sparse optimization iterations and error tolerance. 
Step 5 shows the numerical algorithm of optimization for our model (Eq. \ref{Model2}). Step 5.1 deletes useless the Gaussian basis functions if the corresponding coefficient $c_i$ of ellipsoid Gaussian RBF is less than $1e^{-3}$ per 20 steps. Step 5.4 control the parameter $w_s$ and $w_l$ adaptively by balancing the values of $E_s$ and $E_{l1}$. When $E_s$ is relative small, put more optimization efforts on $E_{l1}$, otherwise focus on minimization of $E_s$. Step 5.5 checks the maximum of error between $\tilde{\phi}$ and $\phi$ at constrained points $\left\{\mathbf{y}_{k}\right\}_{k=1}^{M}$ and corrects the coefficients $w_s$ and $w_{l}$. Step 5.6, after doing $SparseNiter$ iterations, with the number of effective basis is fixed, keep doing some steps of minimization of Es to achieve better accuracy of the approximation on constrained points. Step 5.7-5.8 show the frame of gradient descent method by inexact line search. 


\section{Results and discussion}\label{section4}
In this section, we present some numerical experimental examples to illustrate the effectiveness of our model and method for presenting the Gaussian surface sparsely. Comparisons are made among our model, the original definition of Gaussian molecular surface and sparse RBF method \cite{2016LiSparse}. A set of biomolecules taken from the RCSB Protein Data Bank is chosen as a benchmark set. The number of atoms in these biomolecules ranges from hundreds to thousands. These molecules are chosen randomly from RCSB Protein Data Bank, and no particular structure is specified. All computations were run on a computer with Intel Xeon CPU E5-4650 v2, 2.4GHz, and 126GB memory under a 64-bit Linux system. Further quantitative analysis of the result is given in the following subsections.

\subsection{Spare optimization results}
Fifteen biomolecules are chosen to be sparsely represented by ellipsoid Gaussian RBFs using our model and sparse RBF method \cite{2016LiSparse}.  For fair comparison, the initialize center of RBFs are selected to atom center coordinates for both methods. Table \ref{ProteinResult} shows the final number of effective basis from the results of our method and sparse RBF method.

\begin{table}[H]
\begin{center}
\caption{Number of atoms for 15 test proteins. The third line shows the number of RBFs by sparse RBF method. The last line shows the number of ellipsoid Gaussian RBFs by our method. The results focus on Gaussian molecular surface of decay rate d in Eq. \ref{GaussKernel} equals to 0.5}\label{ProteinResult}
\begin{tabular}{c|cccccc}
\hline\hline
PDBID  & DIALA & ADP  & 2LWC & 3SGS & 1GNA \\ \hline
NATOM  & 20    & 39   & 75   & 94   & 163  \\ \hline
Sparse RBF    & 13   & 8  & 51  & 56 & 108 \\ \hline
\textbf{OUR} & \textbf{4}     & \textbf{7}    & \textbf{11}   & \textbf{17}   & \textbf{28}   \\ \hline
\hline
PDBID  & 1V4Z  & 1BTQ & 6BST & 1MAG & 1BWX \\ \hline
NATOM  & 266   & 307  & 478  & 552  & 643  \\ \hline
Sparse RBF    & 198  & 252 & 316 & 502 & 537 \\ \hline
\textbf{OUR} & \textbf{41}    & \textbf{54}   & \textbf{84}   & \textbf{87}   & \textbf{123}  \\ \hline
\hline
PDBID  & FAS2  & 3SJ4 & 3LOD & 1RMP & AChE \\ \hline
NATOM  & 906   & 1283 & 2315 & 3514 & 8280 \\ \hline
Sparse RBF    & 722  & 953 & 1810 & 2871 & 4438 \\ \hline
\textbf{OUR} & \textbf{142}   & \textbf{233}  & \textbf{530}  & \textbf{701}  & \textbf{1636} \\ \hline\hline
\end{tabular}
\end{center}
\end{table}

\begin{figure}[H]
  \centering
  \includegraphics[scale = 0.6]{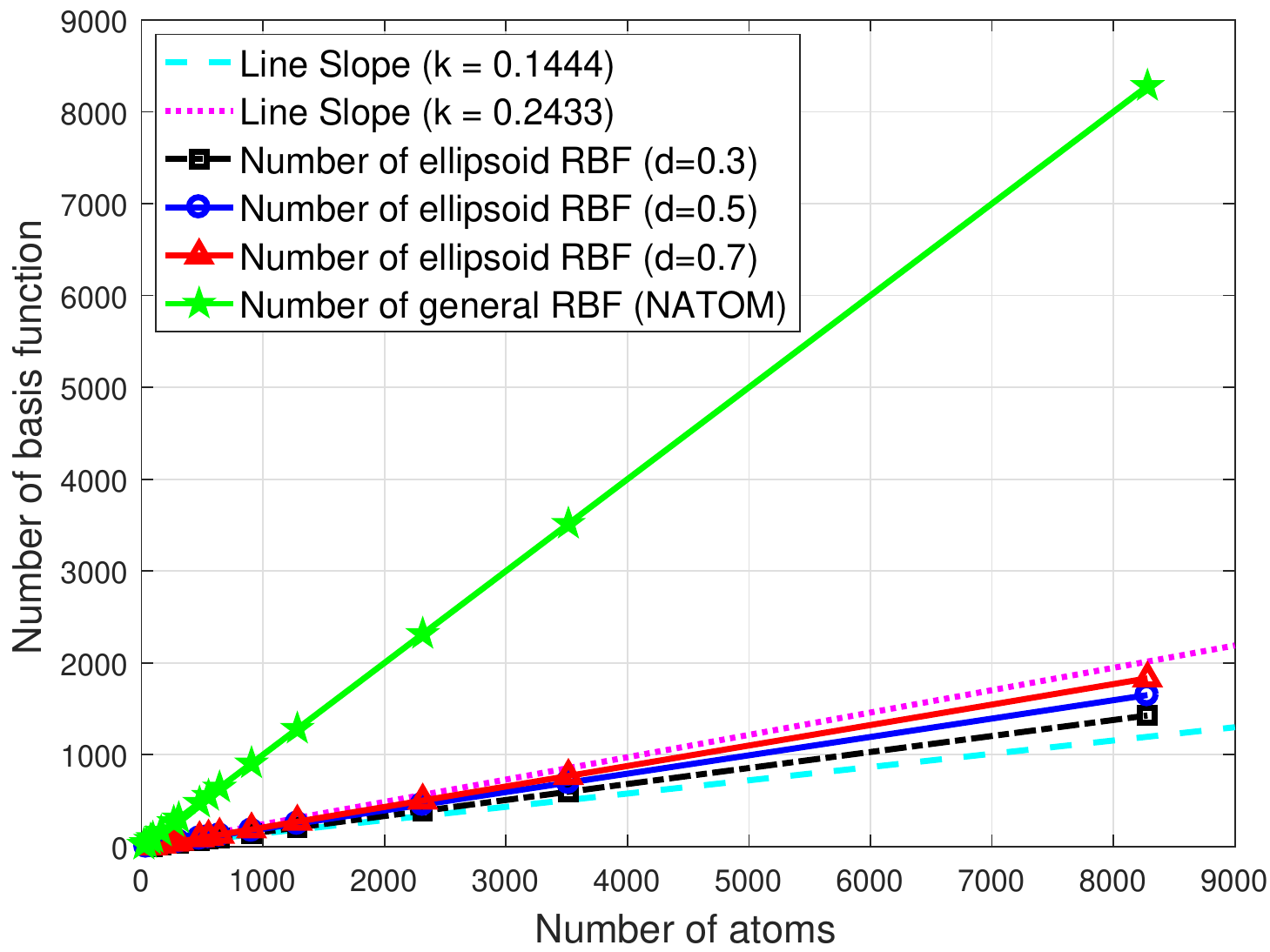}\\
  \caption{Relationship between the number of atoms and the number of ellipsoid Gaussian RBFs after sparse representation.}\label{compressResult}
\end{figure}

Fig. \ref{compressResult} presents the relation between the number of ellipsoid Gaussian RBFs in final sparse representation and the number of atoms in the corresponding molecule. The numbers of general RBFs (number of atoms) for original Gaussian molecular surfaces are shown by green lines with pentagram markers. To present sparse level for our method, we define the sparse ratio $S_r$ as:
\begin{equation}
S_r = \frac{N_{ERBF}}{N_{ATOM}},
\end{equation}
where $N_{ERBF}$ is the number of ellipsoid Gaussian RBFs and $N_{ATOM}$ presents the number of atoms. In Fig. \ref{compressResult}, the changes of sparse ratios with respect to number of atoms for different decay rates ($d$ in Eq. \ref{GaussKernel} equals to 0.3, 0.5 and 0.7) are shown by solid lines with square, circle and triangle markers. The slope of dashed line is the lower bound of sparse ratio ($k= 0.1444$). The slope of dotted line is the upper bound of sparse ratio ($k= 0.2433$). The sparse ratios in the results of our numerical experiments is in ($0.1444,0.2433$). The results show that the larger of decay rate $d$, the more complex is its surface and the smaller is also the sparse ratio. The sparse ratios for Gaussian molecular surface with $d = 0.3$ is smaller than those of Gaussian molecular surface with $d = 0.5$ as shown in Fig. \ref{compressResult}.

\begin{figure}[H]
  \centering
  \includegraphics[scale = 0.6]{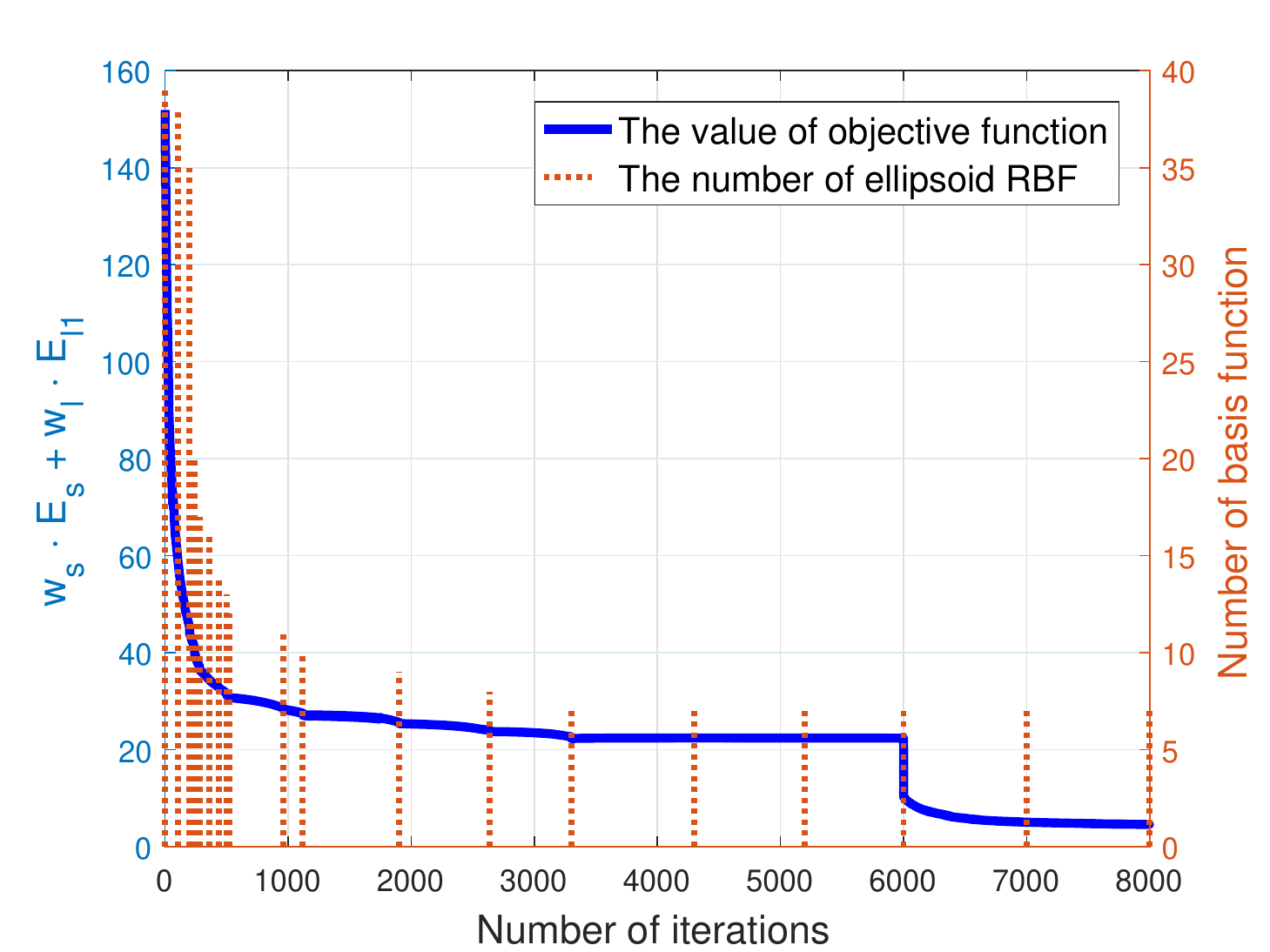}\\
  \caption{One test of the our algorithm on molecule ADP. The blue curve is the objective function trajectory during the 8000 iterations. The red vertical lines represent the number of basis function, whose heights are proportional to the value of the objective function. The number of initial ellipsoid RBF for this trial is 39 and the number of final ellipsoid RBF is 7.}\label{CostNERBF}
\end{figure}

Fig. \ref{CostNERBF} shows the objective function and the number of ellipsoid RBFs is decreasing as the number of iterations increases in the experiment for molecule ADP. In this experiment, the $MaxNiter$ and $SparseNiter$ is set to be 8000 and 6000, respectively. After 6000 iterations, $w_l$ is set to be zero to minimize $E_{l1}$ solely, thus the value of objective function has a abrupt change. The number of ellipsoid RBFs are decreasing dramatically during the iteration process. As shown in Fig.  \ref{CostNERBF}, the model with 7 ellipsoid RBF achieves the lowest errors with a relatively small number of ellipsoid RBFs. Therefore, our method does achieve a desired tradeoff between accuracy and structural sparsity. 

\begin{figure}[H]
  \centering
  \includegraphics[scale = 0.6]{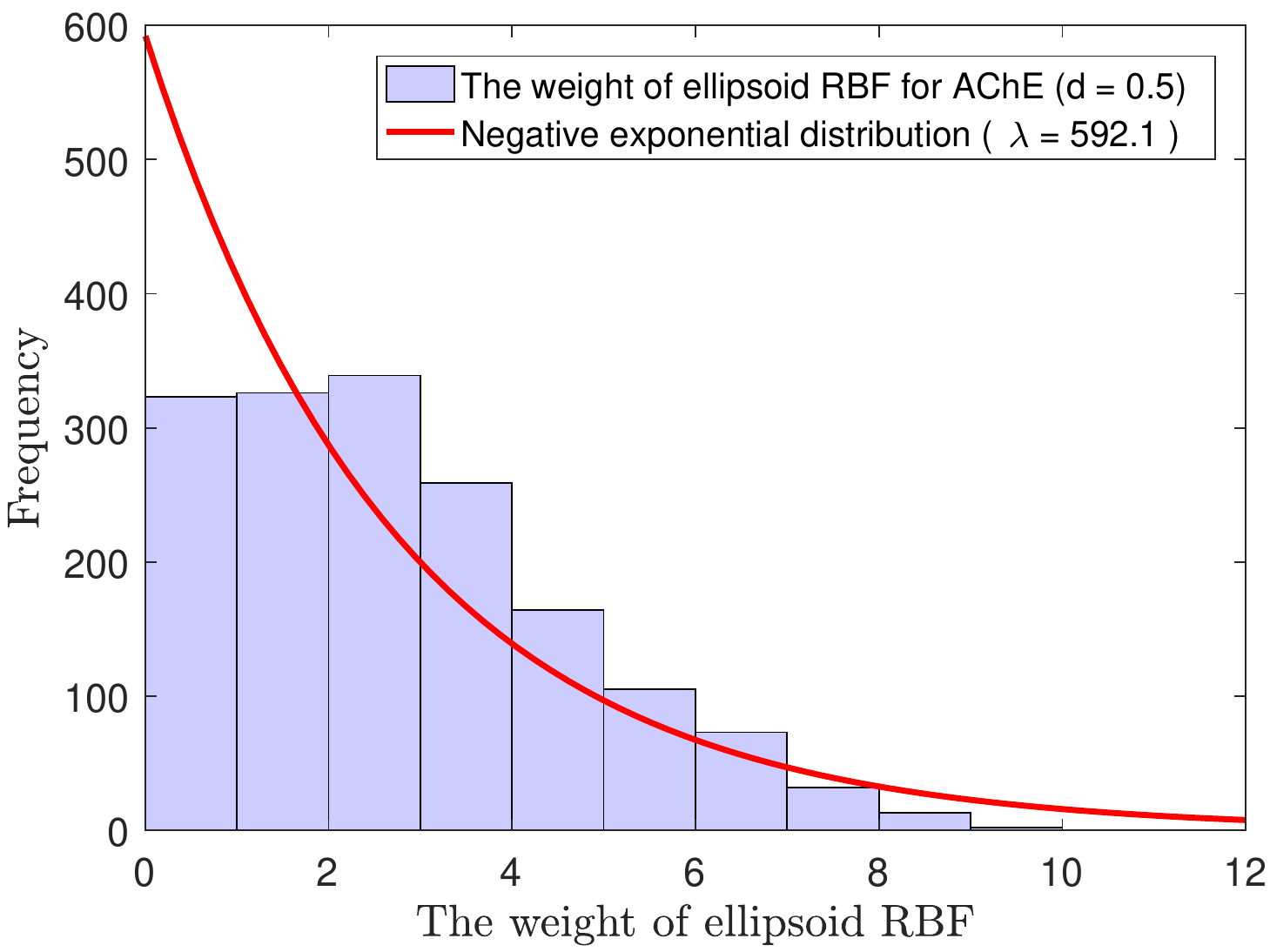}\\
  \caption{Distributions of weight of ellipsoid RBFs in the result of our method for molecule AChE.}\label{HistAChE}
\end{figure}

Fig. \ref{HistAChE} present the weight of ellipsoid RBFs obey negative exponential distribution in the result of our method for molecule AChE. It implies different effectiveness of ellipsoid RBFs after the optimization algorithm.

\subsection{Shape preservation and further results analysis}
In this section, we first check whether the Gaussian surface is preserved after the process of sparse representation through our method. The area of the surface, the volume and the Hausdorff distance are the three criteria to judge whether two surfaces are close enough. These criteria can be calculated on the triangular mesh of the surface. The triangular meshes of comparisons between molecular surfaces before and after sparse representation  are computed through $isosurface$ function in MATLAB. For a triangular surface mesh, the surface area $S$ is determined using the following equation:
\begin{equation}
S=\frac{1}{2} \sum_{i=1}^{n_f}\left\|\overrightarrow {V_{1}^{i} V_{2}^{i}} \times \overrightarrow{V_{1}^{i} V_{3}^{i}}\right\|,
\end{equation}
where $n_f$ is the number of triangle elements and $V_{1}^{i}, V_{2}^{i}, V_{3}^{i}$ denote the coordinates of the three vertices for the $i$th triangle.

The volume V enclosed by the surface mesh is determined using the following equation
\begin{equation}
V=\frac{1}{6} \sum_{i=1}^{n_f} \overrightarrow{V_{2}^{i} V_{1}^{i}} \times \overrightarrow{V_{3}^{i} V_{1}^{i}} \bullet \vec{c}_{i},
\end{equation}
where $c_i$ is the vector from the center of the $i$th triangle to the origin.

The relative errors of area/volume and the Hausdorff distance are used to characterize the difference between the surfaces before and after sparse representation. The relative errors of area and volume are calculated using the following formulas:
\begin{equation}
Error_A = \frac{|A_{our} - A_{original}|}{A_{original}},
\end{equation}
\begin{equation}
Error_V = \frac{|V_{our} - V_{original}|}{V_{original}},
\end{equation}
where $A_{original}$ and $A_{our}$ denote the surface areas of meshes generated from the original and our surfaces respectively. $V_{original}$ and $V_{our}$ denote the corresponding surface volumes of meshes generated from the original and our surfaces respectively.

The Hausdorff distance between two surface meshes is defined as follows.
\begin{equation}
H(S_1,S_2) = \text{max}\left(\mathop{\text{max}}_{p \in S_1}e(p,S_2),\mathop{\text{max}}_{p \in S_2}e(p,S_1)\right),
\end{equation}
where
\begin{equation}
e(p,S) = \mathop{\text{min}}_{p' \in S}d(p,p'),
\end{equation}
S1 and S2 are two piecewise surfaces spanned by the two corresponding meshes, and $d\left(p, p^{\prime}\right)$ is the Euclidean distance between the points $p$ and $p^{\prime}$. In our work, we use Metro \cite{Cignoni1996MME} to compute the Hausdorff distance.

The areas and the volumes enclosed by the surface before and after the sparse representation for each of the molecules are listed in Table \ref{tab:geo}. The Hausdorff distances between the original surface and the final surface for the biomolecules are also listed in Table \ref{tab:geo}.


\begin{table}[H]
\begin{center}
\caption{The areas, volumes and Hausdorff distances obtained with the original and the final surfaces for ten biomolecules. Note: isovalue $\phi = 1.0$, initial decay rate $d = 0.5$.}\label{tab:geo}
\begin{tabular}{|c|c|c|c|c|c|c|c|c|}
\hline
\multirow{2}*{Molecule}  & \multicolumn{3}{c}{Area ($\AA^2$)} & \multicolumn{3}{|c|}{Volume ($\AA^3$)} & \multirow{2}*{Distance ($\AA$)} \\
\cline{2-7}  
 & Original & Our & $Error_A$ & Original & Our & $Error_V$ &   \\ \hline
DIALA & 213.989  & 212.621  & 0.00639 & 252.007  & 250.598  & 0.0056 & 0.287892 \\ \hline
ADP   & 367.853  & 362.241  & 0.01526 & 457.675  & 456.156  & 0.0033 & 0.605255 \\ \hline
2LWC  & 504.97   & 497.506  & 0.01478 & 856.832  & 853.445  & 0.0040 & 0.391414 \\ \hline
1GNA  & 1005.878 & 993.286  & 0.01252 & 1862.528 & 1882.067 & 0.0105 & 0.439894 \\ \hline
1V4Z  & 1479.864 & 1454.448 & 0.01717 & 2834.395 & 2871.164 & 0.0130 & 0.982645 \\ \hline
1BTQ  & 1782.568 & 1757.899 & 0.01384 & 3412.606 & 3471.677 & 0.0173 & 1.134885 \\ \hline
1MAG  & 2479.407 & 2424.142 & 0.02229 & 5730.899 & 5785.473 & 0.0095 & 0.631105 \\ \hline
1BWX  & 2924.699 & 2864.055 & 0.02074 & 6638.106 & 6739.953 & 0.0153 & 0.795473 \\ \hline
FAS2  & 3770.276 & 3664.048 & 0.02818 & 9197.362 & 9301.617 & 0.0113 & 1.292772 \\ \hline
3SJ4  & 5887.047 & 5727.391 & 0.02712 & 13208.53 & 13372.34 & 0.0124 & 0.92122  \\ \hline
\end{tabular}
\end{center}
\end{table}

\begin{figure}[H]
  \centering
  \subfigure[]{
  \label{result:epll}
  \includegraphics[scale = 0.38]{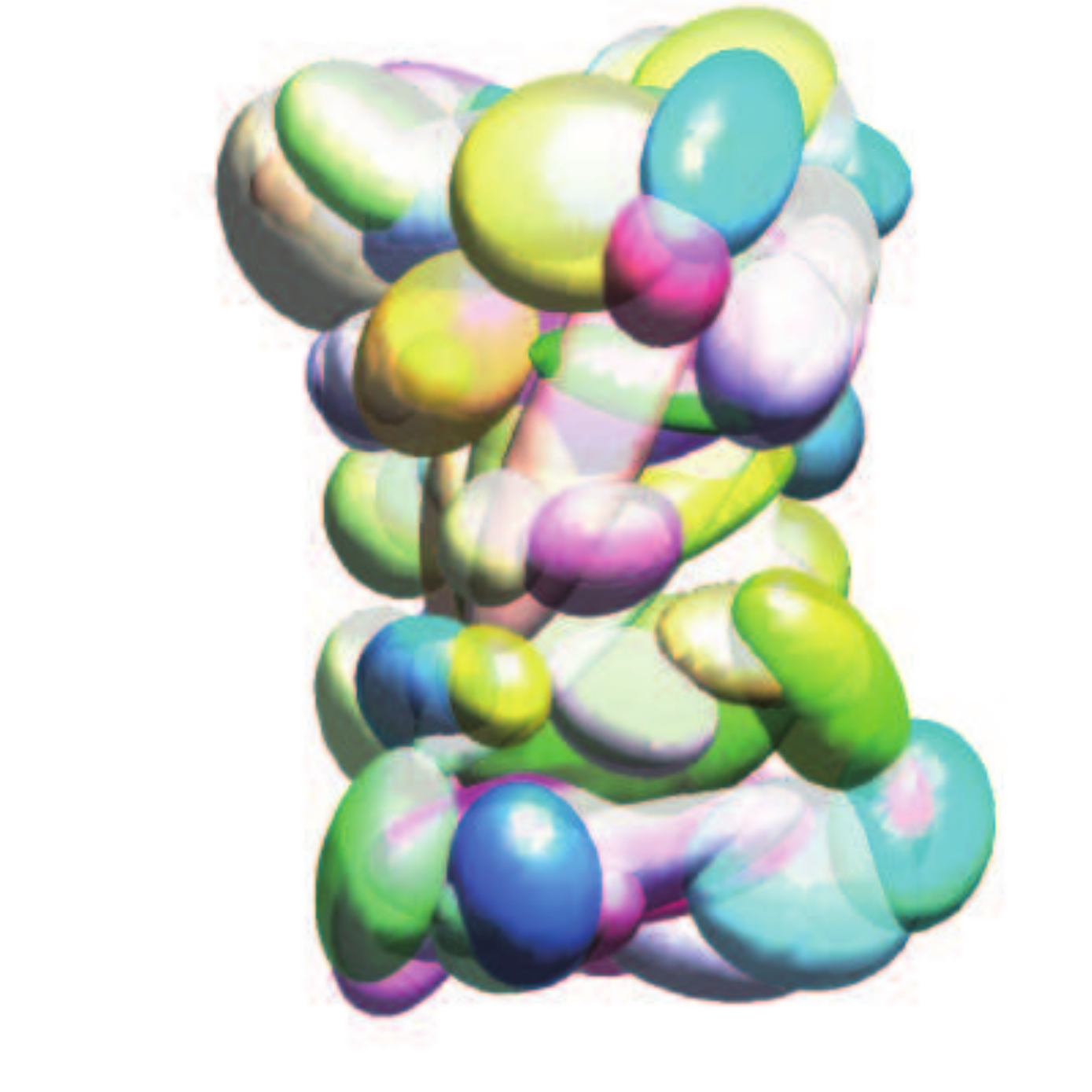}}
  \subfigure[]{
  \label{result:pqr}
  \includegraphics[scale = 0.38]{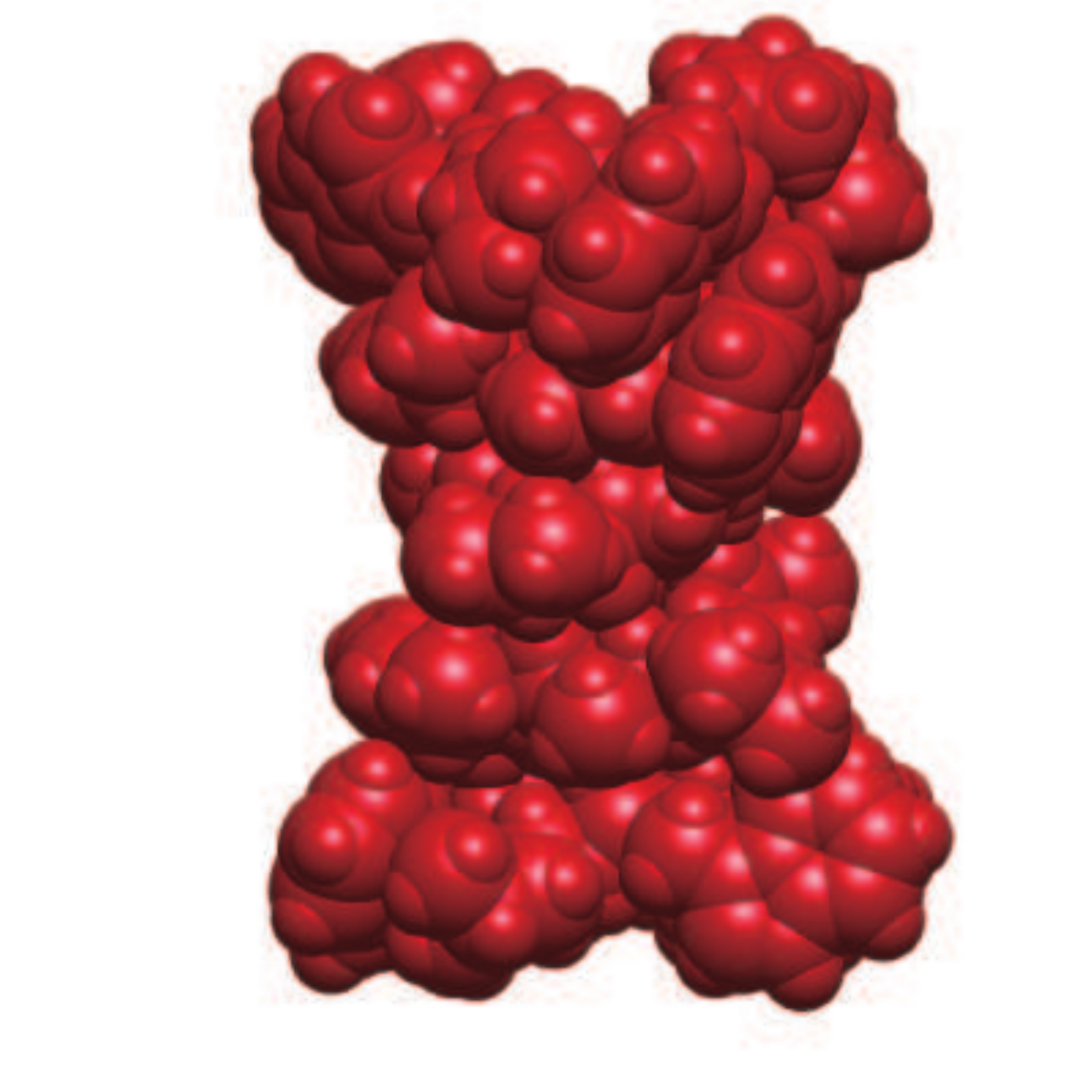}}
  \caption{The locations of the basis functions in our method and the definition of Gaussian surface. This example demonstrates that using our method, the original Gaussian surface is approximated well by a summation of much fewer ellipsoid Gaussian RBFs. (a) The ellipsoid Gaussian RBFs in the sparse representation of molecule 1MAG from our method. (b) The original RBFs in the definition of Gaussian surface of 1MAG. Left: 87 centers, Right: 552 centers.}
  \label{resultCenter}
\end{figure}

Fig. \ref{resultCenter} illustrates the unions of the spheres corresponding to the basis functions by our method and by original surface respectively for the biomolecule (PDBID: 1MAG). This example demonstrates that after the process of sparse representation, the ellipsoid Gaussian RBFs are much sparser than the RBFs in the original definition of Gaussian surface.

\begin{figure}[H]
  \centering
  \subfigure[]{
  \label{geo:volume}
  \includegraphics[scale = 0.6]{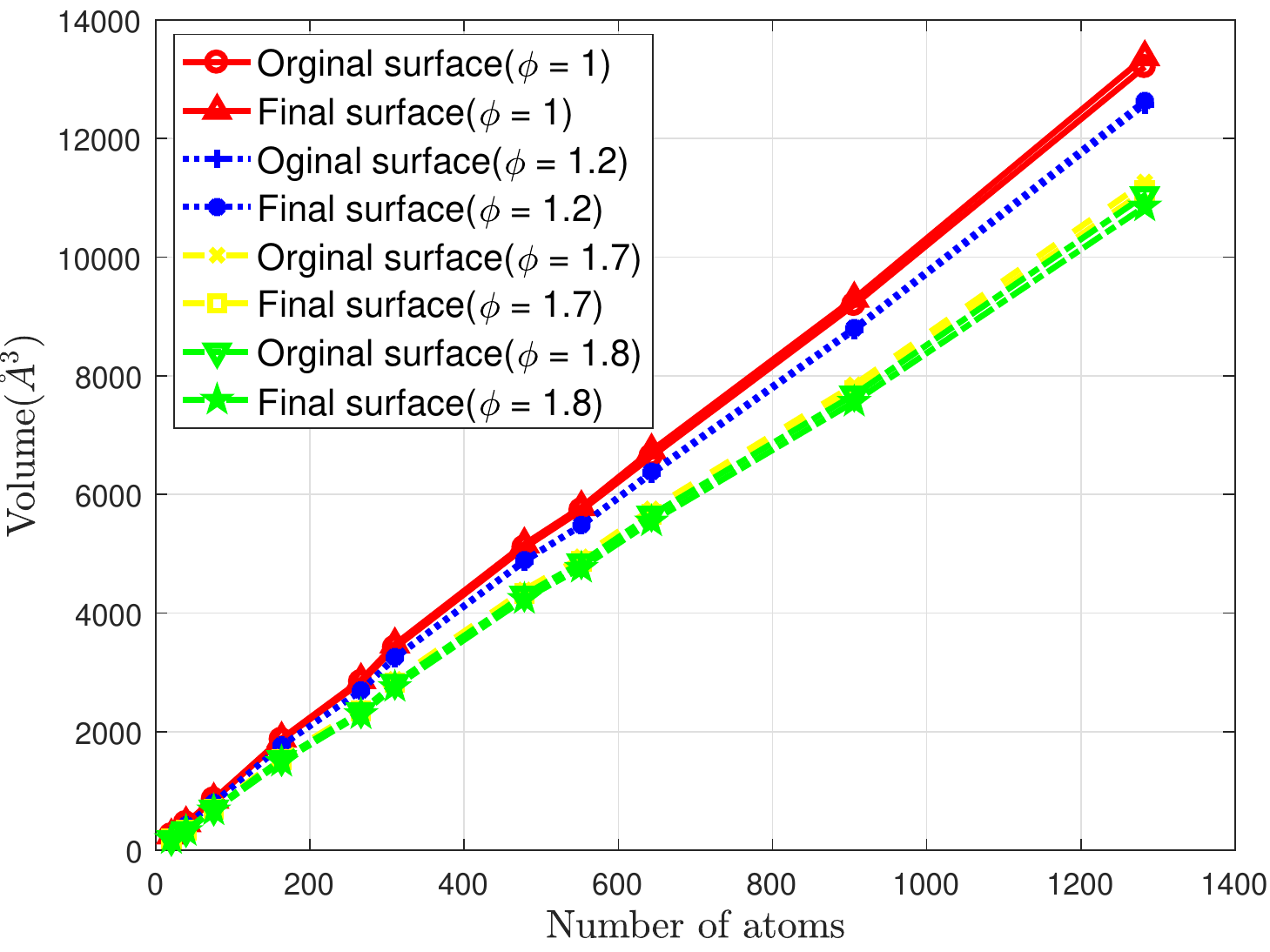}}
  \subfigure[]{
  \label{geo:area}
  \includegraphics[scale = 0.6]{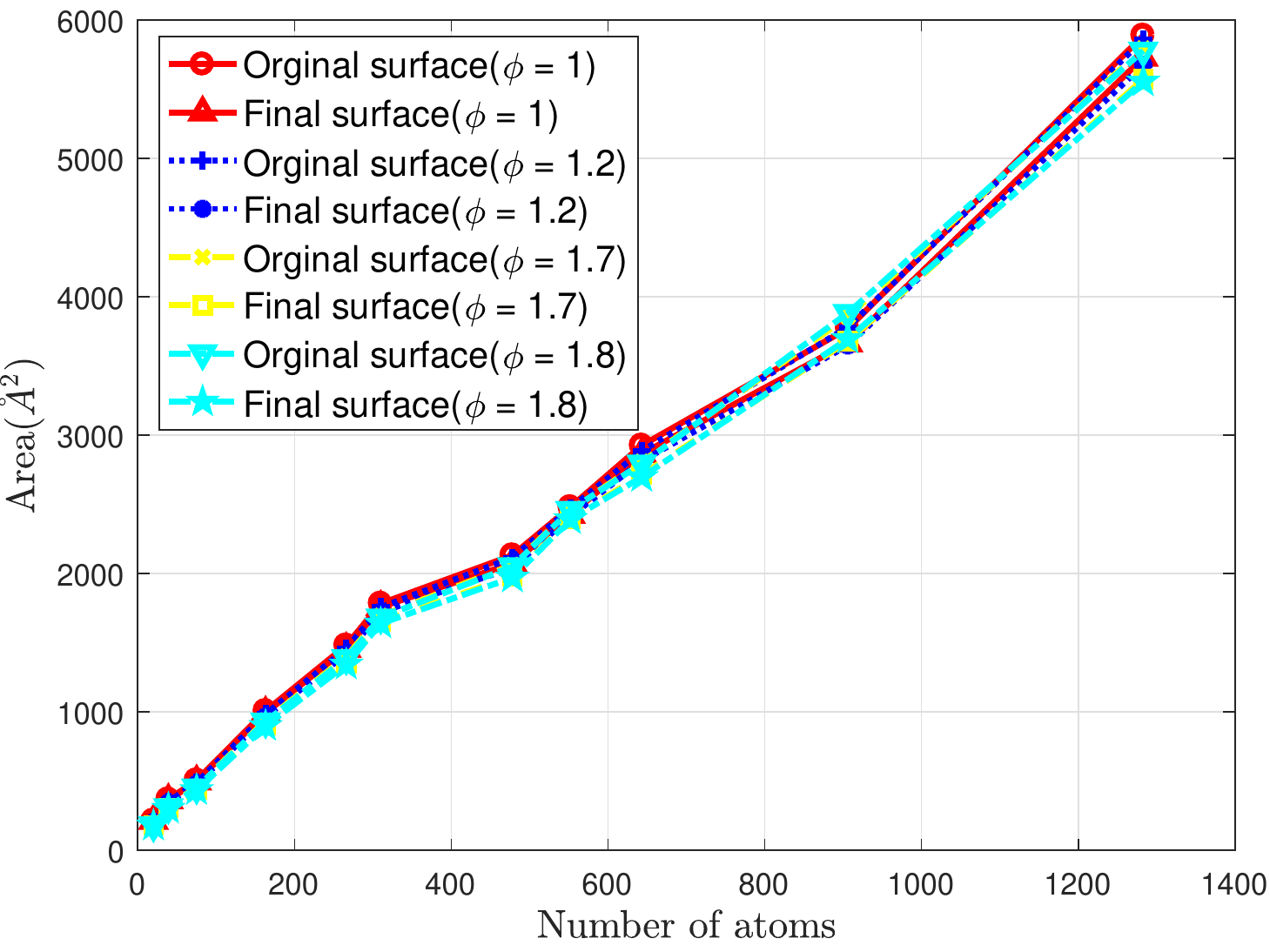}}
  \caption{Volumes and  areas of the initial surface and the sparse surface by our method for the molecules in Table 2. Left: volume; right: area.}
  \label{geo}
\end{figure}

Fig. \ref{geo} plots the areas and volumes of the Gaussian molecular surfaces of different isovalues (c in Eq. \ref{GaussKernel}) before and after spare optimization using our method. And results show that the volumes and areas of Gaussian molecular surfaces with different isovalues are persevered well by our method.

\begin{figure}[h]
  \centering
  \includegraphics[scale = 0.4]{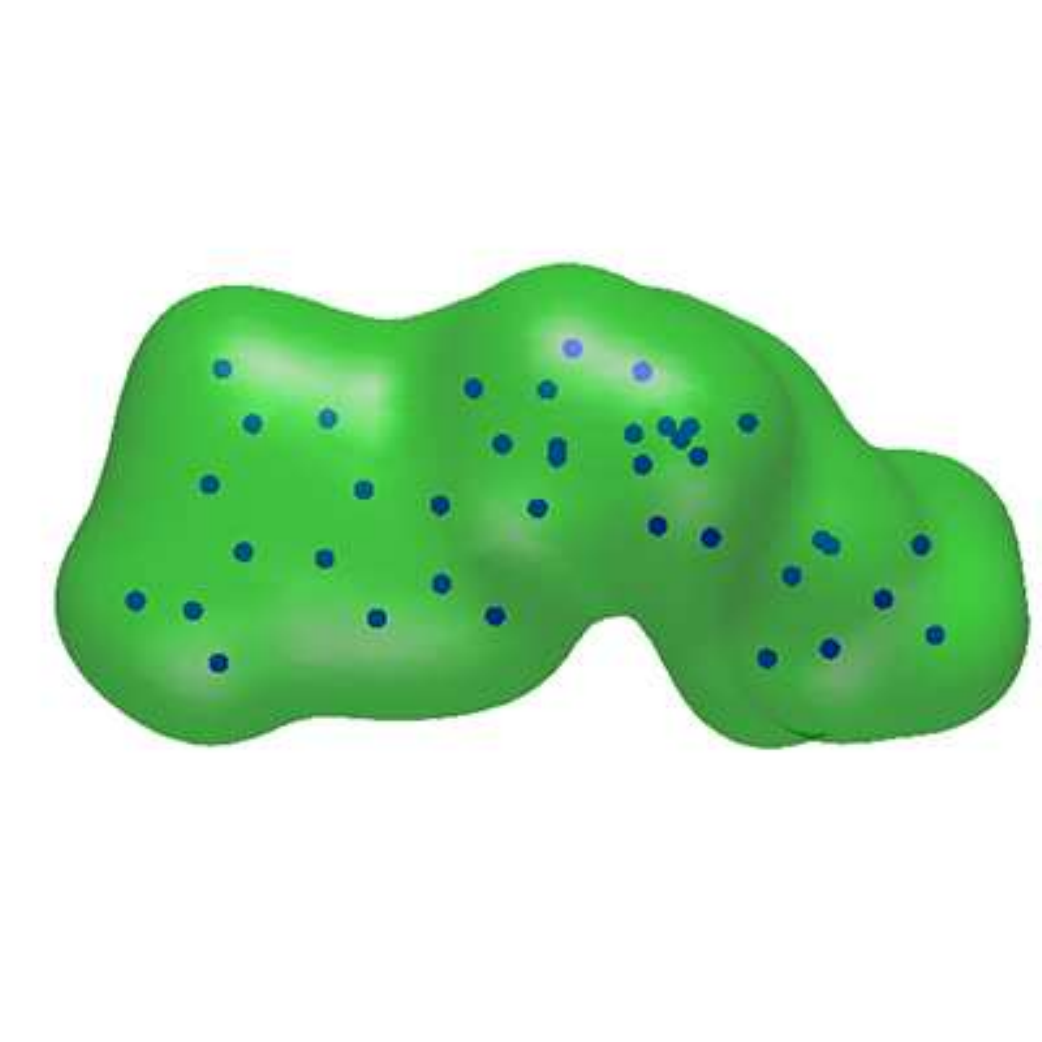}
  \includegraphics[scale = 0.4]{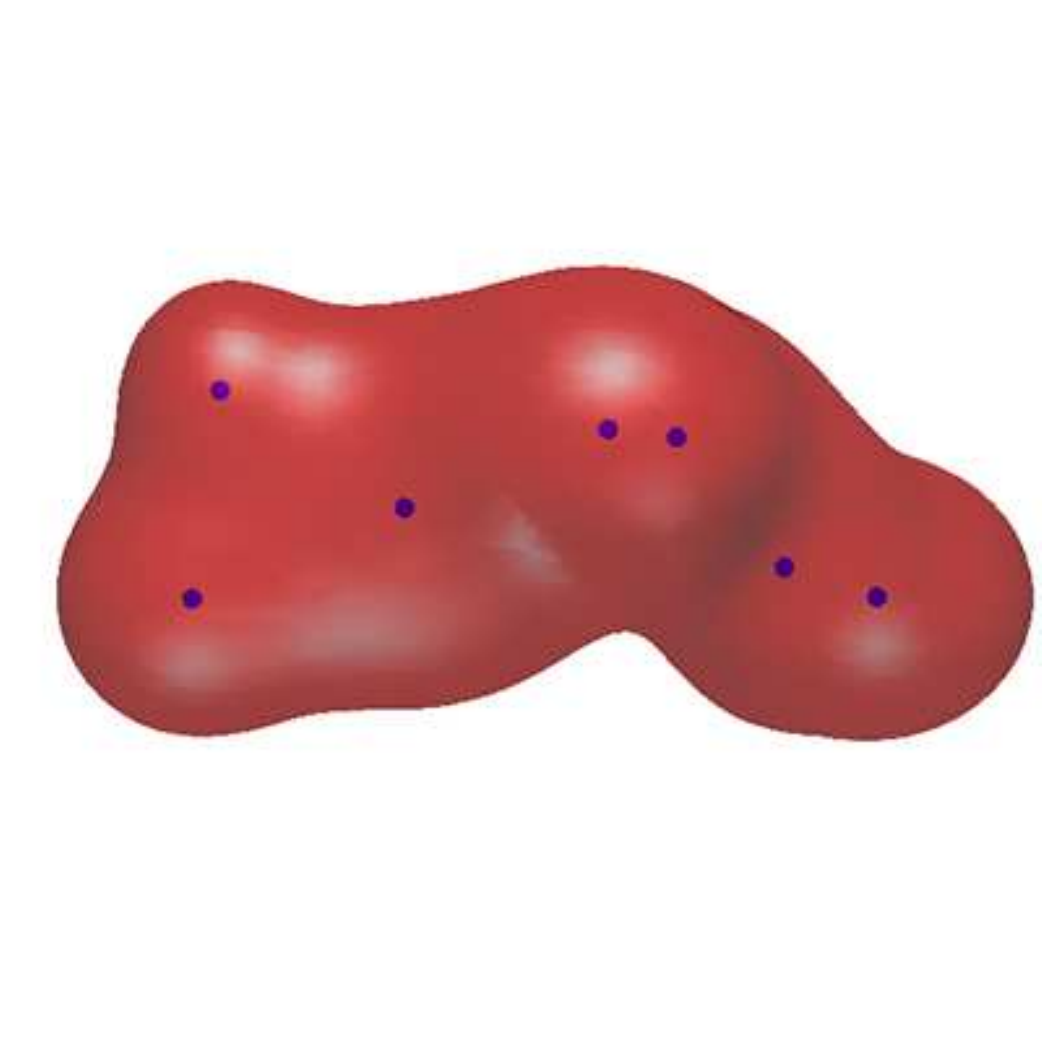}
  \includegraphics[scale = 0.4]{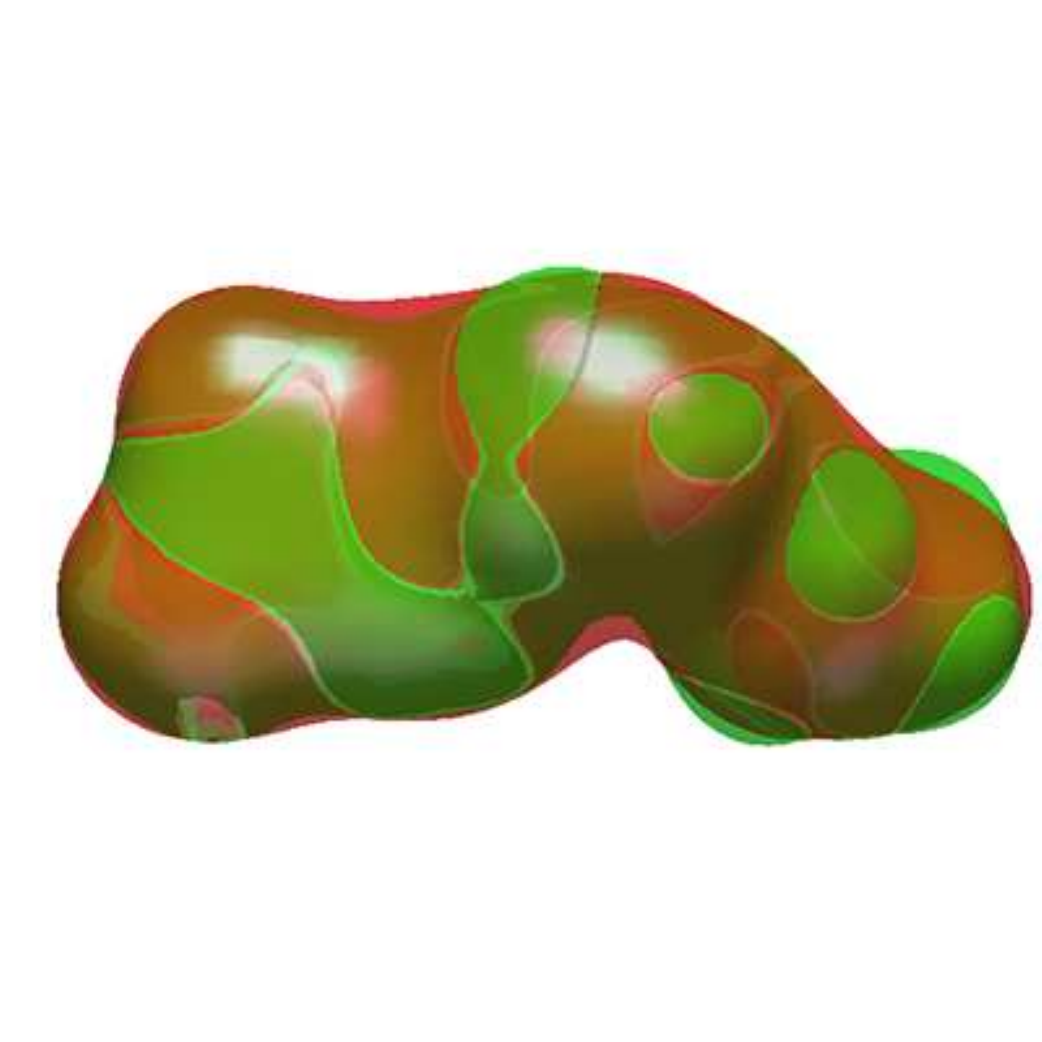}\\
  \includegraphics[scale = 0.4]{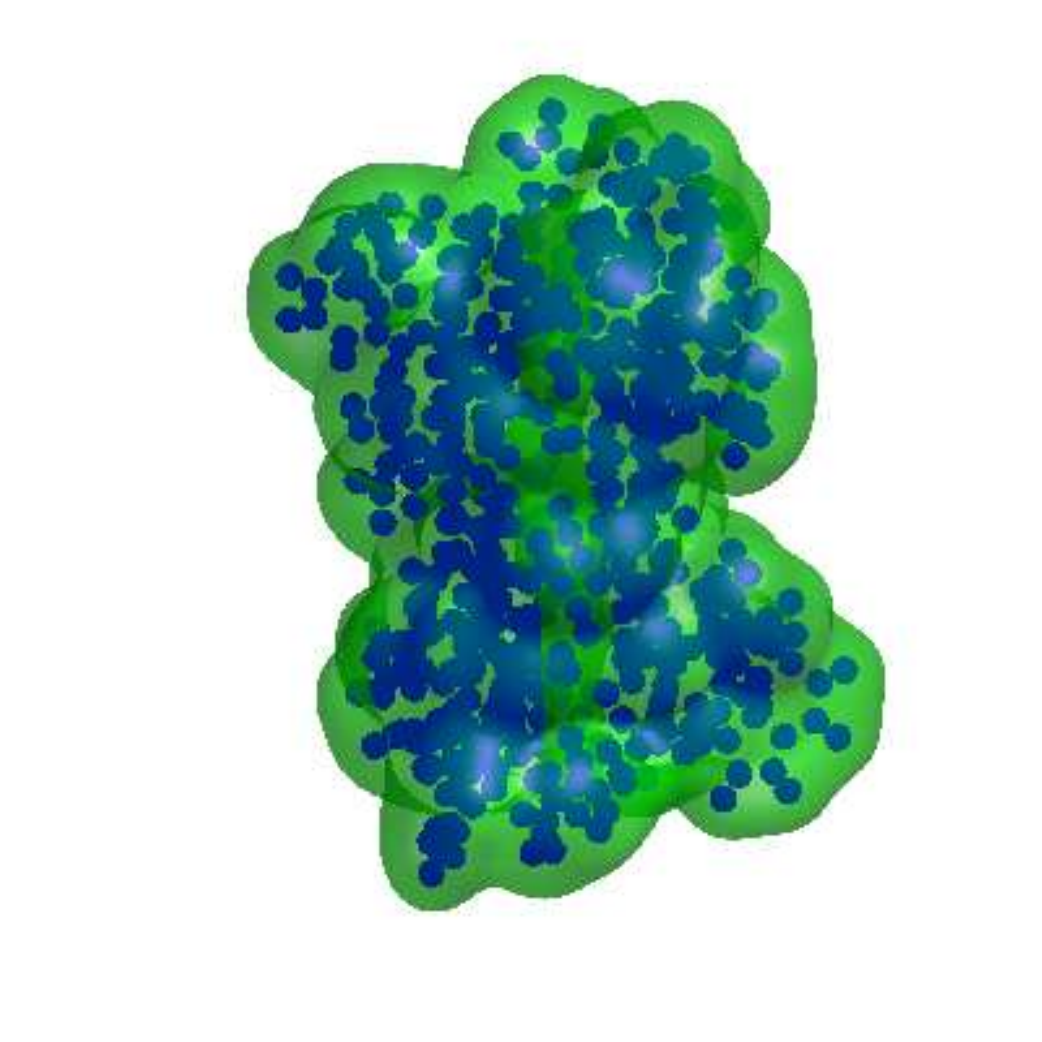}
  \includegraphics[scale = 0.4]{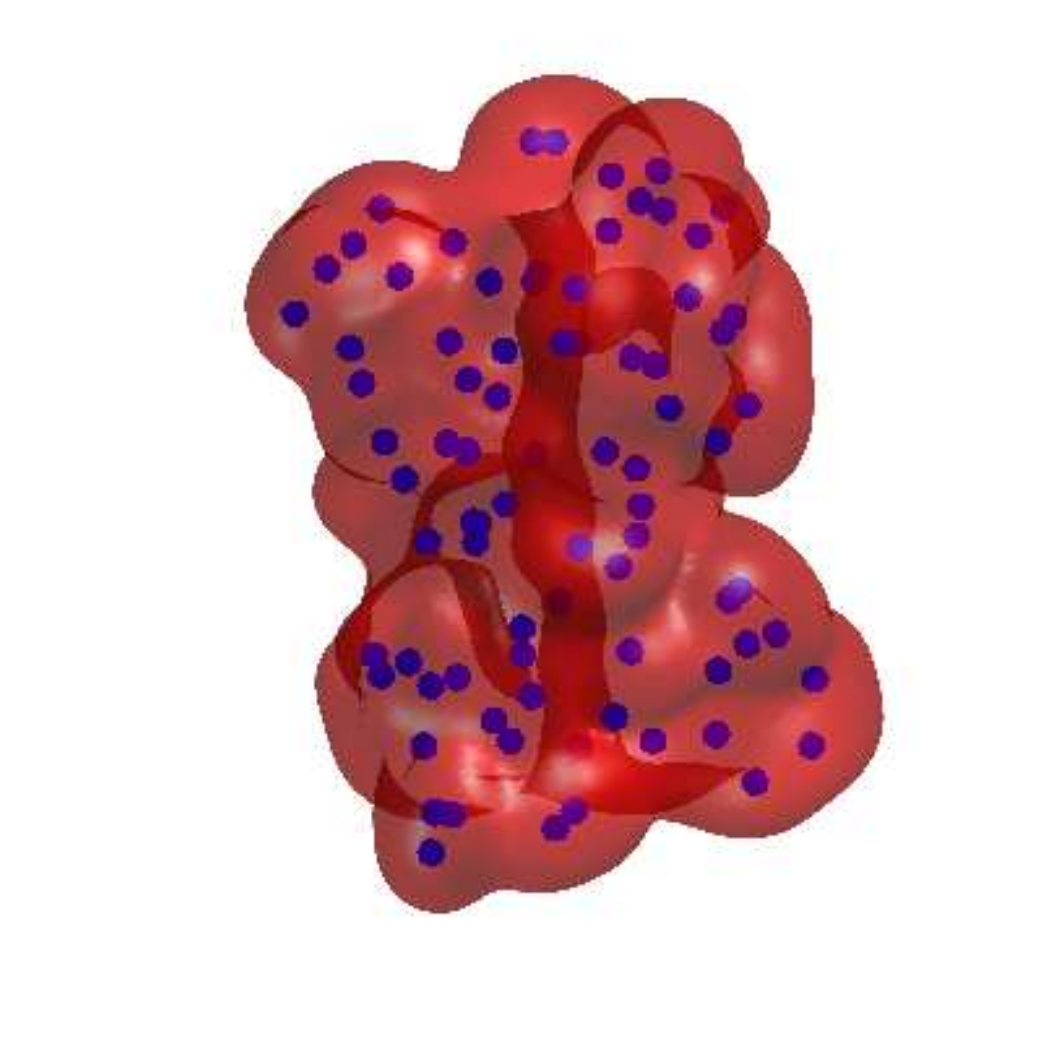}
  \includegraphics[scale = 0.4]{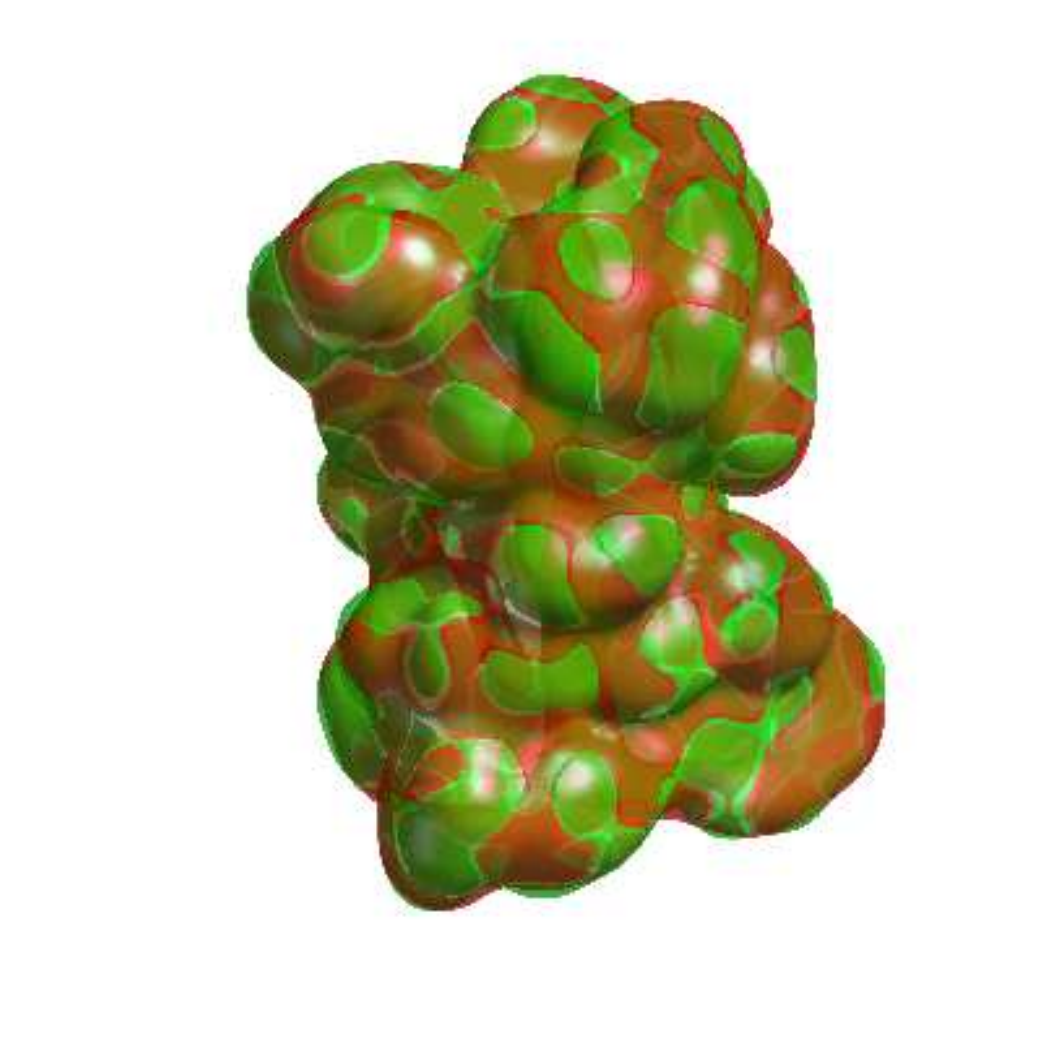}\\
  \includegraphics[scale = 0.4]{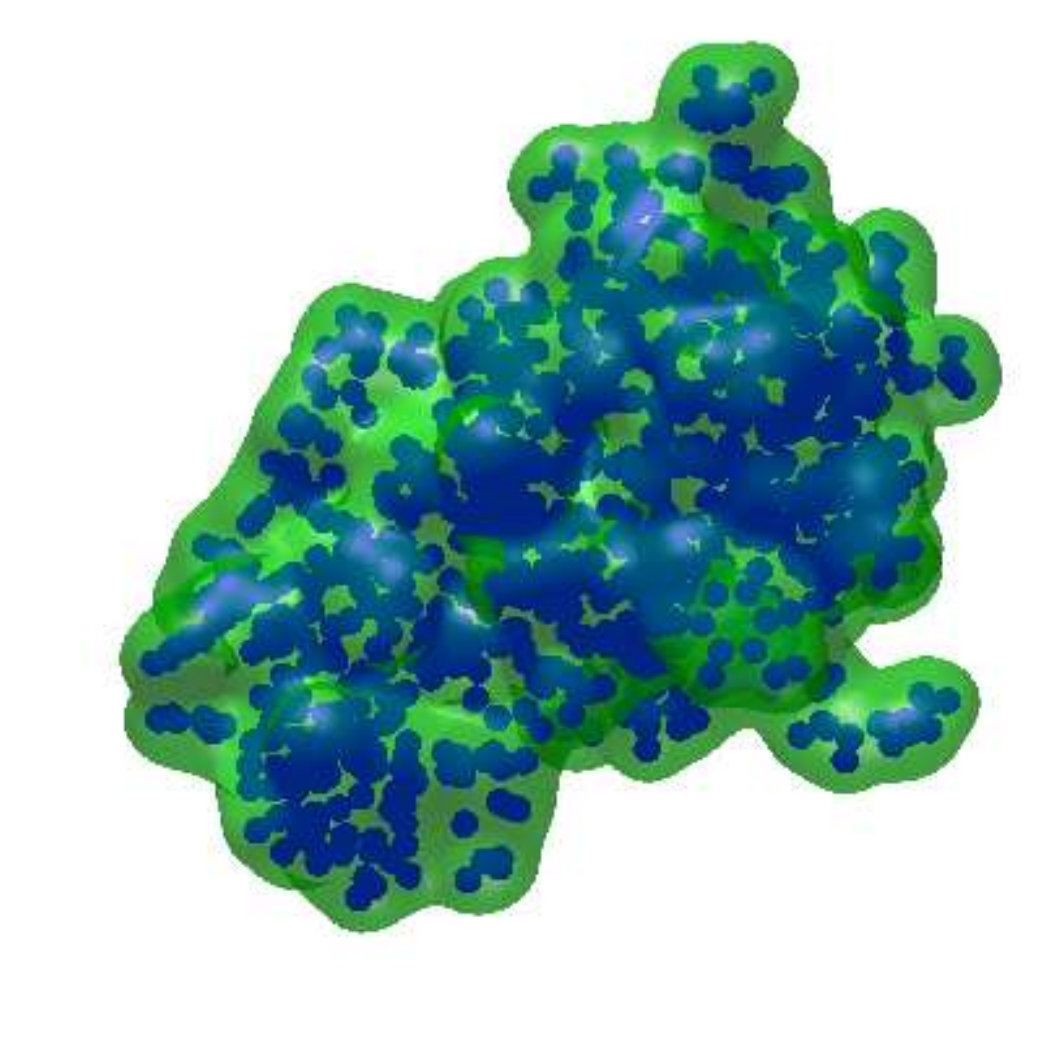}
  \includegraphics[scale = 0.4]{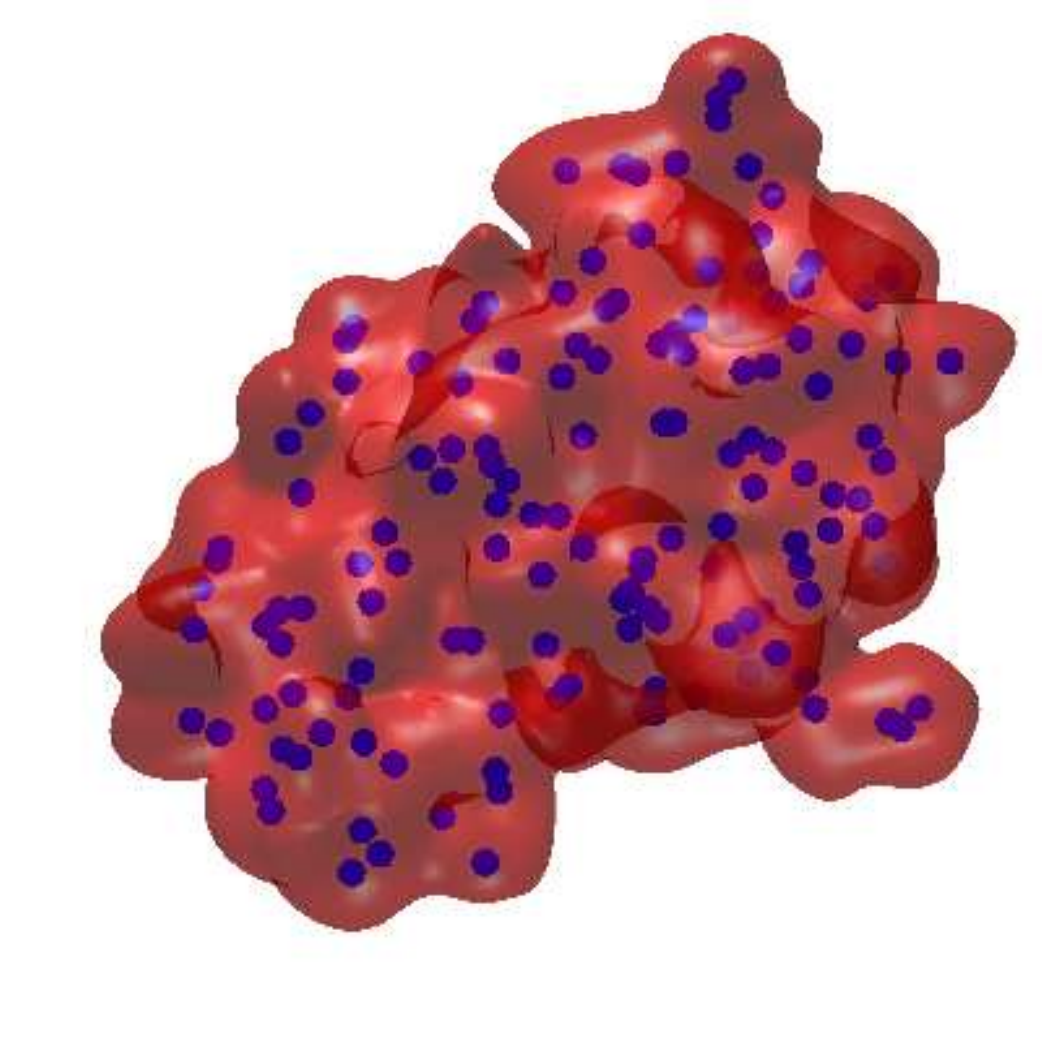}
  \includegraphics[scale = 0.4]{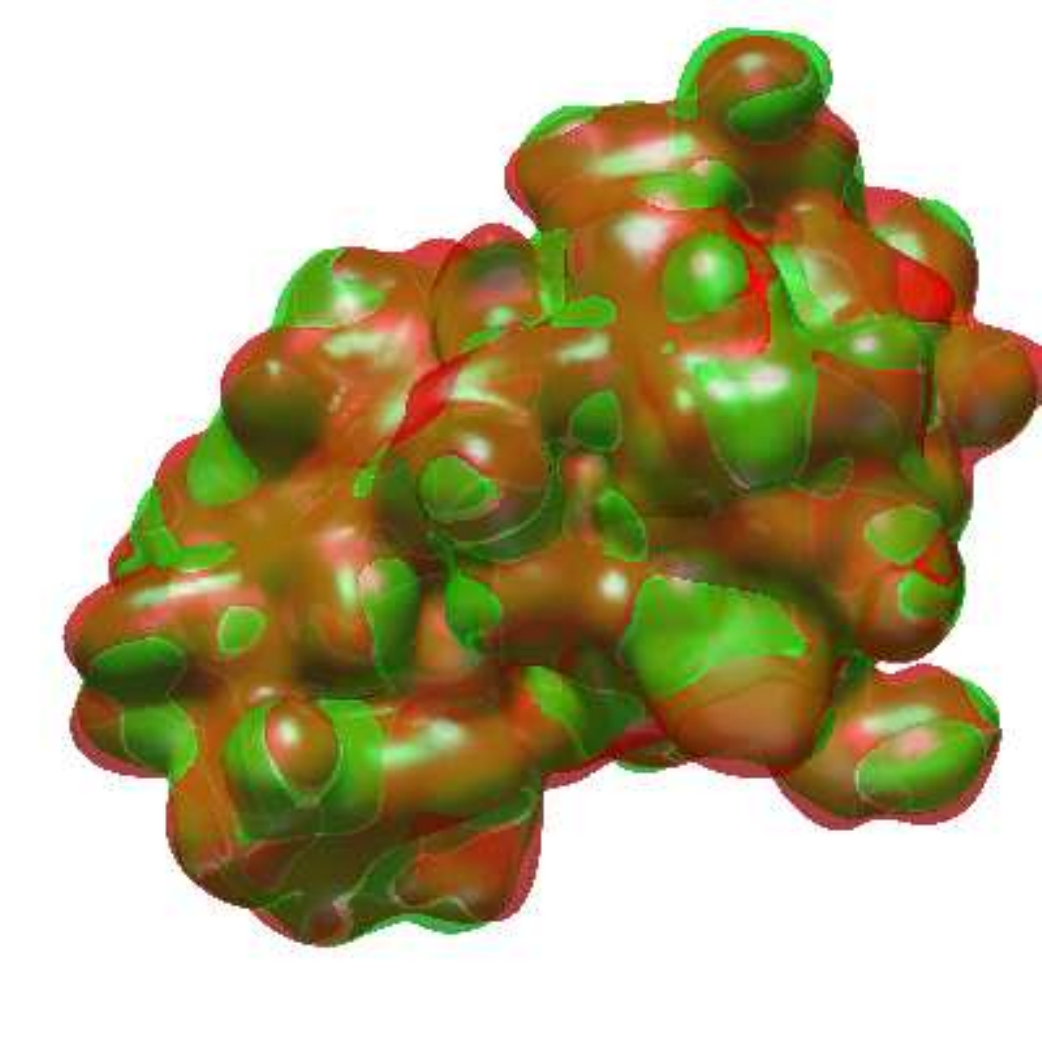}\\
  \caption{Fitting results of our optimization algorithm. Left to right: Original surface (left column), Final surface (middle column) and Original surface overlapped with Final surface(right column). Top to bottom: ADP (first row), 1MAG (second row) and FAS2 (third row). The blue points represent the locations of Gaussian RBF centers.}
  \label{mixresult}
\end{figure}

Fig. \ref{mixresult} illustrates some examples to demonstrate the fitting results of the sparse optimization model. The first column shows original Gaussian surface for three molecules. The second column is the final Gaussian surface in our method, where the blue points represent the location of Gaussian RBF centers. It implies our method need less number of ellipsoid RBFs to represent surface. The last column is the original surface overlapped with the final surface. It presents that the final surface is close to the original surface.

\newpage

\section{Conclusion}\label{section5}
In this paper, a sparse Gaussian molecular surface representation is proposed for arbitrary molecule. The original molecular surface is approximated with ellipsoid Gaussian RBFs. The sparsity of the ellipsoid Gaussian RBF representation is computed by solving an $L_1$ optimization problem. Comparisons and experimental results indicate that our method needs much less number of ellipsoid RBFs to represent the original Gaussian molecular surface.

\section*{Acknowledgements}
The authors thank our group student Qin Wang for her help. The authors thank Manyi Li and Changhe Tu from Shandong university for their help.	
%
%
%
%
%
\bibliographystyle{siamplain}
\bibliography{main}
	
\end{document}